\documentclass[11pt]{article}

\usepackage{multirow}
\usepackage{fullpage}
\usepackage{tabularx}
\usepackage{booktabs}
\usepackage{hyperref}
\usepackage{booktabs}

\usepackage{authblk}

\usepackage[utf8]{inputenc}
\usepackage{subcaption}
\usepackage{nccmath}
\usepackage{graphicx}
\usepackage{multirow}
\usepackage{amsmath,amssymb,amsfonts, amssymb}
\usepackage{mathrsfs}
\usepackage[title]{appendix}
\usepackage{textcomp}
\usepackage{manyfoot}
\usepackage{booktabs}
\usepackage{algorithm}
\usepackage{algorithmicx}
\usepackage{algpseudocode}
\usepackage{listings}
\usepackage{a4wide}
\usepackage[numbers]{natbib}
\usepackage{latexsym}
\usepackage[mathscr]{eucal}
\usepackage{epsfig}
\usepackage{lscape}
\usepackage{amsfonts}
\usepackage{graphicx}
\usepackage{xcolor}
\usepackage{yhmath}
\usepackage{mathdots}
\usepackage{MnSymbol}
\usepackage{epsfig}
\usepackage{lineno}
\usepackage{tikz}
\usepackage{bbold}
\usepackage[absolute]{textpos}

\newtheorem{theorem}{Theorem}
\newtheorem{proposition}[theorem]{Proposition}

\newtheorem{remark}{Remark}

\newcommand{\ZZ}{\mathbb{Z}}
\newcommand{\PP}{\mathbb{P}}
\newcommand{\DD}{\mathbb{D}}
\newcommand{\NN}{\mathbb{N}}
\newcommand{\CC}{\mathbb{C}}
\newcommand{\RR}{\mathbb{R}}
\newcommand{\TT}{\mathbb{T}}

\newcommand{\dedication}[1]{
  \begin{center}
    \textit{#1}
  \end{center}
  \vspace{1em}
}

\title{\textbf{A mixed interpolation-regression method for numerical integration on the unit circle using zeros of para-orthogonal polynomials}}

\author[1]{Ruymán Cruz-Barroso}
\author[2]{Lidia Fernández}
\author[3]{Francisco Marcellán}

\affil[1]{Departamento de Análisis Matemático and Instituto de Matemáticas y sus Aplicaciones (IMAULL), La Laguna University, Spain}
\affil[2]{IMAG and Departamento de Matem\'atica Aplicada, Universidad de Granada, Spain}
\affil[3]{Departamento de Matemáticas, Universidad Carlos III de Madrid, Spain}

\makeatletter
\newcommand{\keywords}[1]{%
  \begingroup
  \par\vspace{0.5em}
  \noindent\textbf{Keywords: }#1\par
  \endgroup
}
\newcommand{\MSC}[1]{%
  \begingroup
  \par\vspace{0.25em}
  \noindent\textbf{MSC (2010): }#1\par
  \endgroup
}
\makeatother

\date{}

\begin{document}

\maketitle

\dedication{To Professor Paul Van Dooren in occasion of his 75th anniversary.}

\begin{abstract}
A new alternative numerical procedure to the  Szeg\H{o} quadrature formulas for the estimation of integrals with respect to a positive Borel measure $\mu$ supported on the unit circle is presented. As in many practical situations, we assume that the values of the integrand $F$ are only known at a finite number of points, which we will assume to be uniformly distributed on the unit circle (although this does not actually constitute a restriction). Our technique consists of obtaining an approximating Laurent polynomial $L$ to $F$ by interpolation in the Hermite sense in a collection of these points that mimic the zeros of a para-orthogonal polynomial with respect to $\mu$, and to use the values of $F$ at the remaining nodes to improve the accuracy of the approximation by a process of simultaneous complex regression. Some numerical examples are carried out.
\end{abstract}


\MSC{33C45; 42C05; 65D32; 41A55; 62J05}

\section{Introduction}\label{sec1}

Let $\omega$ be a positive Borel mesure supported on $[a,b]$, $-\infty \leq a < b \leq \infty,$ and consider the associated inner product
$$\langle f,g\rangle_{\omega}=\int_{a}^b f(x)g(x)d\omega(x), \quad f,g \in L_2^{\omega}=\left\{ f:[a,b] \rightarrow \RR \;: \int_{a}^b \left(f(x)\right)^2d\omega(x) < \infty \right\}.$$
In the estimation of integrals like
\begin{equation}\label{integralGauss}
I_{\omega}(f)=\int_{a}^{b} f(x)d\omega(x),
\end{equation}
quadrature formulas (q.f. from now on) are of the form
\begin{equation}\label{Inomega}
I_n^{\omega}(f)=rA_0f(a)+\left(\sum_{j=1}^{n-r-s}A_jf(x_j)\right) + sA_{n-r}f(b),
\end{equation}
where $r,s \in \{ 0,1 \}$, $x_i \neq x_j$ if $i \neq j$, $x_i \in (a,b)$.

Let us denote by $\PP=\CC[z]$ and $\Lambda=\CC[z,z^{-1}]$ the linear space of polynomials and Laurent polynomials in the variable $z$ with complex coefficients, respectively. $\PP_n$ is the $(n+1)$-dimensional subspace of polynomials of degree less than or equal to $n$ ($\PP_{-1}=\emptyset$) and for two integers $p$ and $q$, $p \leq q$, $\Lambda_{p,q}=\textrm{span}\left\{ z^p,\ldots,z^q \right\}$, whose dimension is $q+1-p$. Let also denote by $\TT=\left\{ z \in \CC : |z|=1 \right\}$ and $\DD=\left\{ z \in \CC : |z|<1 \right\}$ the unit circle and the open unit disk, respectively.

Due to the Weierstrass Approximation Theorem, these q.f. are constructed so that they are exact in a subspace of polynomials of dimension as large as possible. The interpolatory-type rules are obtained from the exact integration of the interpolatory polynomial in Lagrange form at the nodes $\{ x_j \}_{j=1}^{n}$. Thus, with no restriction on the distinct nodes we can obtain in this way exactness for at least $f \in \PP_{n-1}$. The rules (\ref{Inomega}) of maximal degree of accuracy are called of {\em Gauss-type} q.f. ({\em Gaussian}: $r=s=0$; {\em Gauss-Radau}: $r+s=1$; {\em Gauss-Lobatto}: $r=s=1$). They are unique, exact for $f \in \mathbb{P}_{2n-1-r-s}$, positive (with positive weights, that is of importance due to convergence and stability reasons) and the interior nodes are the zeros of the  $(n-r-s)$-th orthogonal polynomial with respect to the positive modified measure $d\omega_{r,s}(x)=(x-a)^r \cdot (b-x)^s \cdot d\omega(x) \geq 0$, $\forall x \in [a,b]$, $r,s \in \{ 0,1 \}$. It is also very well known that an effective computation of those q.f. can be carried out from an eigenvalue problem associated with the corresponding tridiagonal Jacobi matrices, see e.g. \cite{Gau}.

When $\mu$ is a positive Borel mesure supported on $I:=[0,2\pi)$ and the estimation of integrals like
\begin{equation}\label{integralSzego}
I_{\mu}(F)=\int_{0}^{2\pi} F\left( e^{i\theta}\right)d\mu(\theta)
\end{equation}
is required, the analog on the unit circle are q.f of the form
\begin{equation}\label{Inmu}
I_n^{\mu}(F)=\sum_{j=1}^{n} \lambda_jF(z_j), \quad z_j \neq z_k \;\;\textrm{if} \;\; j \neq k, \quad \{ z_j \}_{j=1}^{n} \subset \TT.
\end{equation}

As a consequence of the density of $\Lambda$ in ${\cal C}(\TT)=\left\{ f : \TT \rightarrow \CC \;:\; f \; \textrm{continuous} \right\}$ with respect to the uniform norm (something that does not hold in general for the space $\PP$), it is standard to make use of Laurent polynomials when dealing with approximation problems on the unit circle. The interpolatory-type rules are deduced in this case from the exact integration of the interpolatory Lagrange-Laurent polynomial at the nodes $\{ z_j \}_{j=1}^{n}$, so with no restriction on the distinct nodes we can obtain exactness in subspaces of the form $\Lambda_{p,n-1-p}$, with $0 \leq p \leq n-1$. Szeg\H{o} q.f., introduced for the first time in \cite{JNT} along with invariant para-orthogonal polynomials, are of the form (\ref{Inmu}), of interpolatory-type, exact in a maximal subspace of Laurent polynomials and positive. Since then and  motivated by numerous applications, there are many contributions in the literature concerning numerical integration on the unit circle and properties of invariant para-orthogonal polynomials (that are the analog on the unit circle of orthogonal polynomials associated with positive Borel measures supported on the real line). It is a fact that Szeg\H{o} q.f. are not so well known as Gaussian rules, so let us summarize here the most relevant properties needed for our purposes.

Consider the induced inner product
$$\langle f,g\rangle_{\mu}=\int_{0}^{2\pi} f(z)\overline{g(z)}d\mu(\theta), \quad \quad f,g \in L_2^{\mu}=\left\{ f:\TT \rightarrow \CC : \int_{0}^{2\pi} \left|f(z)\right|^2d\mu(\theta) < \infty \right\}, \quad z=e^{i\theta},$$
and denote by $\{ \rho_n \}_{n\geq 0}$ the family of {\em monic Szeg\H{o} polynomials} with respect to $\mu$, that is, $\rho_n \in \PP_n \backslash \PP_{n-1}$ for all $n\geq 0$ and $\rho_n \perp_{\mu} \PP_{n-1}$. These polynomials satisfy the forward recurrence relation (\cite[Theorem 11.4.2]{Sz} or \cite[Theorem 1.5.2]{SimonBk})
$$\left(\begin{array}{c}
\rho_{n}(z)\\\rho_{n}^{*}(z)\end{array}\right)=\left(\begin{array}{cc}
z&\delta_{n}\\ \overline{\delta_{n}}z&1\end{array}\right)\left(\begin{array}{c}
\rho_{n-1}(z)\\\rho_{n-1}^{*}(z)\end{array}\right),\quad n\geq 1, \quad \rho_{0}=\rho_{0}^{*} \equiv 1,$$
where $\rho_n^{*}(z)=z^n\overline{\rho_n\left( 1 / \bar{z} \right)}\in \PP_n$ are the reciprocal (or reversed) polynomials and
\begin{equation}\label{Verblcoeff}
\delta_0=1, \quad \quad \delta_n=\rho_n(0)=-\frac{\langle z\rho_{n-1}(z),1\rangle_{\mu}}{\langle \rho_{n-1}^{*}(z),1\rangle_{\mu}}\in \DD, \quad \textrm{for} \; n\geq 1,
\end{equation}
are the {\em Verblunsky coefficients}\footnote{There are at least another four different terminologies: Schur, Geronimus, reflection and Szeg\H{o} coefficients, see discussion in \cite[Section 1.1]{SimonBk}.}. It is very well known (see \cite{Sz}) that the zeros of $\rho_n$ lie in $\DD$, so they can not be used, as in the real line case, as nodes of a q.f. on the unit circle.

As it is very well known, polynomials on the real line with real coefficients have their zeros on the real line or appearing in complex conjugate pairs $\{\alpha,\overline{\alpha} \}$. On the unit circle we have a similar property making use of invariant polynomials, that means $B_n^{*}=\tau B_n$ (this definition only has sense for $\tau \in \TT$). This is equivalent to the fact that the zeros of an invariant polynomial are located on the unit circle or they appear in symmetric pairs with respect to $\TT$: $\left\{\alpha, 1/\overline{\alpha} \right\}$. It also follows from definition that an invariant polynomial of exact degree $n$ is orthogonal to $z^k$, if and only if, it is orthogonal to $z^{n-k}$, $1 \leq k \leq n-1$.

Invariant para-orthogonal polynomials with respect to $\mu$ satisfy $B_n \perp_{\mu} \textrm{span} \left\{ z,\ldots,z^{n-1} \right\}$. They were characterized in \cite[Theorem 6.1]{JNT} in terms of the corresponding Szeg\H{o} polynomials:
\begin{equation}\label{para}
B_n(z,\tau_n)=c_n \left[ \rho_n(z) +\tau_n \rho_n^{*}(z) \right], \quad \quad c_n \neq 0, \quad \quad \tau_n \in \TT.
\end{equation}
We will assume $c_n=1$ without loss of generality. Their zeros  are simple and located on $\TT$ (\cite[Theorem 6.2]{JNT}). Moreover, the zeros of $B_n(z,\tau_n)$ and $B_n(z,\tilde{\tau}_n)$ interlace on $\TT$, for $\tau_n,\tilde{\tau}_n \in \TT$, $\tau_n \neq \tilde{\tau}_n$ (see \cite[Theorem 1]{Gol}). These zeros are precisely the nodes in the q.f. (\ref{Inmu}) of maximal degree of accuracy. Indeed, for fixed $\tau_n \in \TT$, positive weights $\lambda_j$, $j=1,\ldots,n, $ can be determined so that
\begin{equation}\label{SantosNjastad}
I_{\mu}(L)=I_n^{\mu}(L), \quad \quad \textrm{for all} \;\; L \in \Lambda_{-(n-1),n-1} \cup \textrm{span}\left\{ \frac{z^{n}}{\delta_{n}+\tau_{n}}-\frac{1}{(\overline{\delta_{n}+\tau_{n}})z^{n}} \right\},
\end{equation}
that is a linear subspace of Laurent polynomials of dimension $2n$ (see \cite[Section 7]{JNT} and \cite{Santos,CDP}). An effective computation of Szeg\H{o} q.f. can be also carried out from an eigenvalue problem, details can be found e.g. in \cite{RuyD}.

In many practical applications the function $f$ in (\ref{integralGauss}) is not known at each point of the interval but only at a finite number of nodes (a discrete grid data), often equispaced. As indicated in \cite{Del4}, the task of reconstructing a function $f$ starting from values obtained by devices, or experimental measures of various kinds arises very often in many different scientific contexts. In such situations, several alternative techniques of numerical integration have been developed in recent years, mainly by F. Dell'Accio, F. Marcellán, F. Nudo and collaborators.

The technique has been presented e.g. in \cite{Del1} in order to avoid the effect of Runge phenomenon when algebraic polynomial interpolation on $N+1$ uniformly distributed nodes ($N \in \NN$) is used by considering mock-Chebyshev interpolation. Suppose $m,r\in \NN$, $0 < m \leq r \leq N$. This procedure generates a polynomial $P \in \PP_r$ that interpolates the function $f$ on a subset of $m+1$ points of these uniformly distributed nodes that mimic the Chebyshev-Lobatto points, and that simultaneously minimizes the approximation error at the discarded points. This procedure was successfully extended to the zeros of orthogonal polynomials instead of the Chebyshev-Lobatto nodes in \cite{Del8}, presenting therefore a much more generic process. Some other extensions have been carried out, for example, by considering Hermite interpolation (see \cite{Del9}) or quasi-uniform grids (see \cite{Del4}). The bivariate case has been also considered by making use essentially the same ideas; see \cite{Del3,Del6,Del7}.

The growing interest in problems on the unit circle in last decades demands the search for alternative numerical integration techniques. The aim of this paper is to introduce for the first time in the literature a mixed interpolation-regression method for numerical integration, in a similar way as it was done in particular in \cite{Del8} for the real line case. Here, the techniques must be modified, mainly due to the fact that we need to work with complex instead of real variables; also the use of Laurent polynomials instead of ordinary polynomials and properties of para-orthogonal polynomials instead of orthogonal polynomials are required. To explore this attractive problem was actually suggested in \cite{CMich}.

It is true that on the unit circle, the Runge phenomenon does not typically occur due to the interpolation by trigonometric functions. The Gibbs phenomenon can take place when dealing with discontinuous functions on the unit circle (which is distinct from the Runge phenomenon, but of a similar nature). In the context of interpolation on the complex plane, in general, a Runge-type phenomenon could occur if the function under consideration is not analytic in a region near the domain, and the convergence of the interpolants may depend on the domain of analyticity of the function. A typical example is the meromorphic function $F(z)=\frac{1}{z-\alpha}$ with $\alpha \in \DD$. Polynomial interpolants based on equally spaced nodes on the unit circle may not converge to $F$ uniformly, and can diverge on points on the unit circle close to the pole $\alpha$ (a Runge-type phenomenon caused by near singularities).

There is another motivation for our problem. Suppose that the values of the function $F$ in (\ref{integralSzego}) are known only in an equispaced set of points on the unit circle. These points correspond to the set of zeros of a para-orthogonal polynomial related to the Lebesgue measure $d\mu(\theta)=d\theta$. So, in this case, interpolatory-type q.f. will produce an effective procedure in the estimation of $I_{\mu}(F)$ since such a q.f. is actually a Szeg\H{o} rule, that is positive and the convergence can be assured for functions $F$ defined on the unit circle that are Riemann integrable. But for general positive Borel measures $\mu$ supported on $[0,2\pi)$, intepolatory-type rules based on these equispaced nodes may not be an effective procedure because the weights in (\ref{Inmu}) are in general complex numbers, and the convergence and stability of the procedure can not be always guaranteed. This motivates clearly to find alternative techniques of numerical integration on the unit circle when the values of the function $F$ in (\ref{integralSzego}) are known only in a set of distinct points on the unit circle, that we will suppose equispaced, although we must point out here that this is not actually a restriction in the development of this work.

The paper has been organized as follows. In Section \ref{sec2} we present a mixed interpolation-regression method for the numerical estimation of integrals of the form (\ref{integralSzego}). The idea is to construct an approximating function $L$ to $F$ in a linear subspace of Laurent polynomials of dimension $r+1$ that interpolates $F$ in a collection of $m$ points that mimic the zeros of a para-orthogonal polynomial of degree $m$ with respect to the measure $\mu$, and that uses some of the discarded nodes to improve the accuracy of the approximation by a process of simultaneous complex regression. We explain in Section \ref{sec3} how to construct the interpolating part, consisting in an interpolating Laurent polynomial in Hermite sense, where the number of derivatives that the function $L$ matches at each point is arbitrary. Some illustrative numerical examples and conclusions are shown in Sections \ref{sec4} and \ref{sec5}, respectively.

\section{A mixed interpolation-regression method}\label{sec2}

Suppose that the parameters  $m,r,N\in \NN$, $0 < m \leq r \leq N,$ are fixed in advance. Let $Z_{N}=\left\{ z_j \right\}_{j=0}^{N-1} \subset \TT$ be an ordered set of $N$ uniformly distributed nodes on the unit circle $\TT$ labelled counterclockwise:
\begin{equation}\label{ZN}
z_j=e^{i\left( \theta_0+hj \right)}, \quad \quad h=\frac{2\pi}{N}, \quad \quad j=0,\ldots,N-1, \quad \quad \textrm{for some fixed } \theta_0 \in I=[0,2\pi).
\end{equation}
As it was said before, we could actually work with an arbitrary distribution of $N$ distinct nodes on the unit circle, but we will assume a uniform grid for the shake of simplicity.

Let $F$ be an unknown function defined on $\TT$. We assume to know only its evaluations on the set of points $Z_N$; our purpose is to estimate $I_{\mu}(F)$ given by (\ref{integralSzego}) from the only data available of the function $F$, $\mu$ being a positive Borel measure supported on $I$.

Since the zeros of $\rho_m$ lie in $\DD$, the function $F_m: \TT \rightarrow \TT$ given by $$F_m(z)=-\frac{\rho_m(z)}{\rho_m^{*}(z)}$$ is a finite Blaschke product of degree $m$. Define $\tau_{m}=F_m(z_0)\in \TT$, consider the para-orthogonal polynomial $B_{m}(z,\tau_m)$ given by (\ref{para}), and let us denote its $m$ distinct zeros by $\Xi_{m}=\left\{ \xi_k \right\}_{k=0}^{m-1} \subset \TT$, labelled counterclockwise such that $z_0=\xi_0$.

We assume that it is possible to identify a set $\Upsilon_m=\left\{ \tilde{z}_k \right\}_{k=0}^{m-1} \subset Z_{N}$ such that $\tilde{z}_0=z_0=\xi_0$, $\tilde{z}_i \neq \tilde{z}_j$ for all $i\neq j$, $i,j\in \{ 0,\ldots,m-1 \},$ and that for each $k=1,\ldots,m-1$, the point $\tilde{z}_k$ is a solution of the problem
\begin{equation}\label{problemk}
\displaystyle{\min_{j=0,\ldots ,m-1} \left| \xi_j-\tilde{z}_k \right|, \quad \quad \forall k=1,\ldots,m-1.}
\end{equation}
For a fixed $k$, this problem has always a unique solution unless accidentally the zero $\tilde{z}_k$ of the para-orthogonal polynomial $B_{m}(z,\tau_m)$ is placed just in the middle of the arc in between two consecutive nodes of the uniform grid, namely $z_s,z_{s+1}$. In such a case, we can arbitrarily choose one of these two points to take part of the set $\Upsilon_m$. What is important here is to be sure that a set $\Upsilon_m$ with $m$ distinct points can be constructed with $m-1$ of the solutions of the $m-1$ minimizing problems (\ref{problemk}) along with the point $\tilde{z}_0$. Thus, the points of the set $\Upsilon_m$ are some of the original uniformly distributed nodes chosen as closest to the zeros of $B_{m}(z,\tau_m)$, and we can also order them counterclockwise.

The parameter $m$ will be computed thus as the greatest value so that this property holds. In the real line situation this problem was studied first for the Chebyshev-Lobatto nodes (see \cite{BXu}) by making use of the explicit formulas for such nodes, obtaining thus an explicit expression for $m=m(N)$. In the more general situation considered in \cite{Del8}, where the nodes are zeros of orthogonal polynomials, the question was solved just by noticing that the zeros of the orthogonal polynomials commonly used in applications, such as Legendre and Jacobi, share the same density as the zeros of Chebyshev polynomials of the first kind. However, it was remarked there that for positive Borel measures on the real line having a distribution differing from the Chebyshev points, no general formula for such a parameter $m$ can be obtained, and this observation is not only true for the case of the real axis but also for the unit circle. Some numerical experiments concerning a weight function on the unit circle that depends on a parameter $q\in (0,1)$ have been carried out in Section \ref{sec4} (see further). We have observed that the behavior of $m=m(N,q)$ depends entirely on the parameter $q$ since a radical change in the distribution of the zeros of para-orthogonal polynomials occurs when $q$ is varied.

Our purpose is to approximate $F$ by a Laurent polynomial $L$ that belongs to a linear subspace  $\Lambda_{-p,r-p-1}$, of dimension $r$, where $p$ is a parameter fixed in advance and $0 \leq p \leq r-1$. We will work in general, but a reasonable choice for the parameter $p$ is to take $p=l,$ when $r=2l+1,$ and either $p=l-1$ or $p=l,$ when $r=2l$ (balanced or CMV ordering). The approximating function $L$ to the function $F$ will be constructed in such a way that it interpolates $F$ in the set $\Upsilon_m$ and the fact that the original information that has not been used (that is, the nodes in $Z_N \backslash \Upsilon_m$) will be employed to improve the accuracy of the approximation by a process of simultaneous complex regression. The question of how to determine appropriately the dimension $r$ of the space of Laurent polynomials where the approximating function $L$ belongs will be considered at the end of this section.

Thus, for $m<r$, $L$ is given by
\begin{equation}\label{L}
L(z)=P_m(z)+\frac{\omega_m(z)}{z^p} \cdot Q_{N-m}(z), \quad \quad P_m(z) ,\; \frac{\omega_m(z)}{z^p} \cdot Q_{N-m}(z) \in \Lambda_{-p,r-p-1}.
\end{equation}
Here $P_m$ is an interpolatory Laurent polynomial to the function $F$ at the set of nodes $\Upsilon_m$,
$$\omega_m(z)=\prod_{k=0}^{m-1} \left(z-\tilde{z}_k \right)\in \PP_m$$ is the nodal polynomial and $Q_{N-m}$ an algebraic polynomial
\begin{equation}\label{QN-m}
Q_{N-m}(z)=\sum_{l=0}^{r-m-1}C_l\upsilon_l(z) \in \PP_{r-m-1},
\end{equation}
$\left\{ \upsilon_l \right\}_{l=0}^{r-m-1}$ being some basis of $\PP_{r-m-1}$. For our purposes we will take $\upsilon_l(z)=z^l$, $l=0,\ldots,r-m-1$. Notice that these are the analog on the unit circle of the Chebyshev polynomials defined on $[-1,1]$, since  $T_n(\cos \theta)=\cos (n\theta)=\Re (z^n)$ ($T_n$ being the $n$-th Chebyshev polynomial of the first kind), that is $$T_n \left( \frac{z+z^{-1}}{2} \right)=\frac{z^n+z^{-n}}{2}, \quad \quad z \in \TT.$$

An appropriate election of $P_m$ will be discussed in Section \ref{sec3}. The $r-m$ constants $C_l$ in (\ref{QN-m}) are determined by solving
\begin{equation}\label{QN-mProblem}
\displaystyle{\min_{L\in \Lambda_{p,r-p-1}^{*}} \parallel F-L \parallel_2^2},
\end{equation}
where $\parallel \cdot \parallel_2$ is the discrete 2-norm on $Z_N \backslash \Upsilon_m=\left\{ \hat{z}_s \right\}_{s=1}^{N-m}$ and
$$\Lambda_{p,r-p-1}^{*}=\left\{ L\in \Lambda_{p,r-p-1} \;:\; L\; \textrm{interpolates the function } F \textrm{ in } \Upsilon_m \right\}.$$

\begin{theorem}\label{QN-mProblemThm}
The problem \eqref{QN-mProblem} has a unique solution.
\end{theorem}

\begin{DT}
Observe that (\ref{QN-mProblem}) reads
$$
\begin{array}{ccl}
\displaystyle{\min_{L\in \Lambda_{p,r-p-1}^{*}} \parallel F-L \parallel_2^2}&=&\displaystyle{\min_{Q_{N-m}\in \PP_{r-m-1}} \sum_{s=1}^{N-m} \left| F(\hat{z}_s)-\left( P_m(\hat{z}_s)+\frac{\omega_m(\hat{z}_s)}{\hat{z}_s^p} \cdot Q_{N-m}(\hat{z}_s)\right) \right|^2} \\ \\
&=&\displaystyle{\min_{Q_{N-m}\in \PP_{r-m-1}} \sum_{s=1}^{N-m} \left| \left\{ \frac{\left[F(\hat{z}_s)- P_m(\hat{z}_s)\right]\cdot\hat{z}_s^p}{\omega_m(\hat{z}_s)} - Q_{N-m}(\hat{z}_s)\right\}\cdot \frac{\omega_m(\hat{z}_s)}{\hat{z}_s^p} \right|^2} \\ \\
&=&\displaystyle{\min_{Q_{N-m}\in \PP_{r-m-1}} \parallel G - Q_{N-m} \parallel_{2,\omega_m}^2,}
\end{array}$$
where $\parallel \cdot \parallel_{2,\omega_m}$ is the weighted discrete 2-norm on $Z_N \backslash \Upsilon_m$,
$$\parallel \varphi \parallel_{2,\omega_m}=\left(\sum_{s=1}^{N-m} \left| \varphi(\hat{z}_s)\omega_m(\hat{z}_s) \right|^2 \right)^{\frac{1}{2}}$$
and
$$G(z)=\frac{\left[F(z)- P_m(z)\right]\cdot z^p}{\omega_m(z)},$$
that is well defined on $Z_N \backslash \Upsilon_m$. The result follows since the original problem (\ref{QN-mProblem}) has been transformed in a classical complex linear least squares problem, that has a unique solution (see e.g. \cite[Sections 1-3]{Miller}).
\begin{flushright}
$\Box$
\end{flushright}
\end{DT}

An explicit formula for the polynomial $Q_{N-m}$ can be obtained. The discrete weighted norm $\parallel \cdot \parallel_{2,\omega_m}$ is induced by the discrete inner product
$$\langle \varphi , \psi \rangle_{2,\omega_m}=\sum_{s=1}^{N-m} \varphi(\hat{z}_s)\overline{\psi(\hat{z}_s)}\left|\omega_m(\hat{z}_s) \right|^2.$$
If we denote by $\left\{ p_n \right\}_{n\geq 0}$ the corresponding family of orthonormal polynomials, then it follows that
$$Q_{N-m} (z)=\sum_{s=0}^{r-m-1} \alpha_s p_s(z), \quad \textrm{with Fourier coefficients } \alpha_s=\langle G , p_s \rangle_{2,\omega_m}.$$

However, we can alternatively proceed by taking into account that the solution can by obtained from the {\em normal equations} related to a linear least squares problem. In our complex context, the same formula is used as in the real situation, just by interchanging transpose by Hermitian (conjugate) transpose (see \cite[Sections 1-3]{Miller}). More precisely,
$${\cal C}=\left({\cal H}^{\dag} {\cal H}\right)^{-1}{\cal H}^{\dag}{\cal G}, $$
where
$${\cal C}=\left( \begin{array}{c} C_0 \\ \vdots \\ C_{r-m-1} \end{array} \right), \quad \quad {\cal G}=\left( \begin{array}{c} G(\hat{z}_1)  \\ \vdots \\ G(\hat{z}_{N-m}) \end{array} \right), \quad \quad {\cal H}=\left( \begin{array}{ccc} \upsilon_0(\hat{z}_1) &\cdots & \upsilon_{r-m-1}(\hat{z}_1) \\ \vdots &&\vdots \\ \upsilon_0(\hat{z}_{N-m}) &\cdots & \upsilon_{r-m-1}(\hat{z}_{N-m})\end{array} \right).$$
Here ${\cal H}^{\dag}=\overline{{\cal H}}^T$ is the Hermitian transpose of ${\cal H}$ and ${\cal H}^{\dag} {\cal H}$ is a nonsingular square matrix.

The following result establishes that in the discrete $2$-norm $\parallel \cdot \parallel_2$, the error obtained in the approximation of the function $F$ using the Laurent polynomial $L$ is smaller than or equal to the one obtained from the interpolatory Laurent polynomial $P_m$.

\begin{theorem}\label{QN-merror}
In the discrete 2-norm $\parallel F-L \parallel_2 \leq \parallel F-P_m \parallel_2$ holds.
\end{theorem}

\begin{DT}
Notice that for all $s=0,\ldots,N-m, $ we get
$$\langle G-Q_{N-m},p_s\rangle_{2,\omega_n}=\langle G,p_s\rangle_{2,\omega_n}-\langle \sum_{s=0}^{r-m-1} \alpha_sp_s,p_s\rangle=0,$$ which implies that $G-Q_{N-m} \perp_{2,\omega_n} Q_{N-m}$, and from the Pythagorean Theorem in a pre-Hilbert space,
$$\begin{array}{ccl}
\parallel F-L \parallel_{2}^2  &=& \parallel G-Q_{N-m} \parallel_{2,\omega_n}^2 = \parallel G \parallel_{2,\omega_n}^2 - \parallel Q_{N-m} \parallel_{2,\omega_n}^2 = \parallel F-P_m \parallel_{2}^2 - \parallel Q_{N-m} \parallel_{2,\omega_n}^2 \\
&\leq&\parallel F-P_m \parallel_{2}^2.
\end{array}$$
Occasionally, it may happen that $\parallel Q_{N-m} \parallel_{2,\omega_n}^2=0$, for example, when either $m=r=N$ or if the interpolatory Laurent polynomial $P_m$ at the nodes $\Upsilon_m$ also matches the function $F$ at the discarded nodes $Z_N \backslash \Upsilon_m$. Thus, the inequality is not strict, in general.
\begin{flushright}
$\Box$
\end{flushright}
\end{DT}

The limit case $m=N$ (and hence, equal to $r$) corresponds to a situation where there exists a bijection between the zeros of the para-orthogonal polynomial $B_{N}(z,\tau_N)$ and the set of nodes $Z_N$. In this case, $Q_{N-m}(z) \equiv 0$ in (\ref{L}). We present some numerical examples of this extreme case in Section \ref{sec4}.

The case $m=r<N$ also implies $Q_{N-m}(z) \equiv 0$, and this case really wouldn't make sense. This means that the discarded nodes $Z_N \backslash \Upsilon_m$ have not any role in the construction of $L$. Since the idea is to use precisely such an information of the function $F$ that has not been used yet in order to get a better approximating function to $F$ than $P_m$, we conclude this section by determining appropriately the dimension $r$ of the linear space of Laurent polynomials where $L$ belongs.

In the real line situation this question has been solved (see \cite{Del1}) by making use of a result in  \cite{Rei} and the explicit expression for the Chebyshev-Lobatto nodes. The same value $r$ has been used for example in \cite{Del8} by making use of the same argument that many of the families of orthogonal polynomials that are commonly used in applications have the same density distribution as the zeros of the Chebyshev polynomials of the first kind. But this is not true in general, and the same situation occurs for the unit circle concerning the zero distribution of para-orthogonal polynomials (see e.g. the results presented in Section \ref{sec4} further).

The argument in \cite{Rei} essentially establishes that, for the unit circle case, the degree $r-m-1$ of the approximating polynomial $Q_{N-m}$ in the least squares nodes $Z_N \backslash \Upsilon_m$ has to be taken as the greatest value such that $r-m \in \{ 3,\ldots,N-m \}$ of these nodes are close to $r-m$ Fejér points, that on the unit circle are precisely $r-m$ uniformly distributed nodes. Here we understand that to consider if a partition is close or not to the uniform one, we need to have at least three nodes. Thus, if we consider the partition $$\Delta \;:\quad 0 \leq \hat{z}_1=e^{i\theta_1} < \hat{z}_2=e^{i\theta_2} < \cdots < \hat{z}_{N-m}=e^{i\theta_{N-m}} < 2\pi,$$
then the problem reduces to find a subpartition of $\Delta$ of cardinality the greatest possible value $r-m$ such that these $N-m$ angles are closer to being uniformly distributed. We describe next a possible simple strategy to address this issue that for small values of $N-m$ is a robust and efficient algorithm. For large values of $N-m$ other approaches using combinatorial optimization techniques, dynamic programming, greedy algorithms, linear clustering or robust regression should be considered.

There is a total of
$$2^{N-m}-\left( \substack{N-m \\ 2 } \right)-\left( \substack{N-m \\ 1 } \right)-\left( \substack{N-m \\ 0 } \right) = 2^{N-m}-\frac{(N-m)(N-m-1)}{2}-(N-m)-1$$ subpartitions of $\Delta$ with at least three nodes. For each of these subpartitions
$$\tilde{\Delta}_k \;:\quad 0 \leq \hat{z}_{t_1}=e^{i\theta_{t_1}} < \hat{z}_{t_2}=e^{i\theta_{t_2}} < \cdots < \hat{z}_{t_k}=e^{i\theta_{t_{k}}} < 2\pi,$$
we can obtain a measure of the deviation of $\tilde{\Delta}_k$ from the uniform one if we compute $K_{\tilde{\Delta}_k} \geq 1$ satisfying
\begin{equation}\label{Kopt}
\displaystyle{\frac{\parallel \tilde{\Delta}_k \parallel}{\left| \theta_{t_{j+1}}-\theta_{t_{j}} \right|} \leq K_{\tilde{\Delta}_k}, \quad \quad \forall j=0,1,\ldots,k-1, \quad \quad \parallel \tilde{\Delta}_k \parallel=\max_{j=0,\ldots,k-1} \left| \theta_{t_{j+1}}-\theta_{t_j} \right|.}
\end{equation}
The optimal choice will therefore be the partition $\tilde{\Delta}_k$ for which the value $K_{\tilde{\Delta}_k}$ is minimum.

\section{Laurent polynomial of general Hermite interpolation on the unit cirle}\label{sec3}

We know that the Laurent polynomial $P_m \in \Lambda_{-p,r-p-1}$ (of dimension $r$) must interpolate the function $F$ in (\ref{integralSzego}) at the set of $m$ points $\Upsilon_m$, where $m<r$. We are interested in the CMV case: $p=l$ if $r=2l+1$ and $p=l-1$ or $p=l$ when $r=2l$. The only data that is available is the evaluation of the function $F$ in the set $\Upsilon_m$, so we can proceed by considering an interpolating Lagrange-Laurent polynomial, namely
\begin{equation}\label{Lagrange-Laurent}
P_m(z)=\sum_{k=0}^{m-1} F(\tilde{z}_k)\ell_k(z)\in \Lambda_{-\tilde{p},m-\tilde{p}-1}, \quad \quad \ell_k(z)=\frac{\omega_m(z)}{(z-\tilde{z}_k)\omega_m^{'}(\tilde{z}_k)}\cdot \left( \frac{\tilde{z}_k}{z} \right)^{\tilde{p}}.
\end{equation}
In this case $P_m$ belongs to a linear  subspace of Laurent polynomials of dimension $m.$ So we can take $\tilde{p}=\tilde{l},$ if $m=2\tilde{l}+1,$ and either $\tilde{p}=\tilde{l}-1$ or $\tilde{p}=\tilde{l},$ when $m=2\tilde{l}$ ($p \geq \tilde{p}$). The functions $\ell_k$ are known as the fundamental Lagrange-Laurent polynomials. Notice that $$\Lambda_{-\tilde{p},m-\tilde{p}-1} \subseteq \Lambda_{-p,r-p-1}.$$

In many applications it is often feasible to know the evaluations of the function $F$ in the set $\Upsilon_m$ along with possible some derivatives at these nodes. In such cases one can improve the accuracy of the process by making use of Hermite instead of Lagrange interpolation. There are some contributions in the literature in recent years on this topic. For example, the authors in \cite{Alicia2} describe the behavior of Hermite interpolants for piecewise continuous functions on the unit circle, analyzing the corresponding Gibbs phenomenon near the discontinuities. Particular Hermite interpolation problems on the unit circle are studied in \cite{Alicia1} and \cite{Alicia3} by considering equally spaced nodal points  and zeros of para-orthogonal polynomials, respectively. The well known Hermite-Fejér Theorem concerning the uniform convergence of Hermite-Fejér interpolants to continuous functions on equally spaced nodes was extended also to the unit circle in \cite{Daruis}.

Let us suppose thus that the values of some derivatives of the function $F$ in the set $\Upsilon_m$ are known. Taking into account that the linear subspace of Laurent polynomials where $P_m$ belongs can be extended from $\Lambda_{-\tilde{p},m-\tilde{p}-1}$ to $\Lambda_{-p,r-p-1}$ (that is, from dimension $m$ to $r$), we propose the following general interpolation problem not considered before in the literature that extends to the Laurent case the situation for algebraic polynomials.

\begin{theorem}\label{TheoHermite1}
For an arbitrary set $\Upsilon_m=\left\{ \tilde{z}_k \right\}_{k=0}^{m-1} \subset \TT$, $\tilde{z}_i \neq \tilde{z}_j$ if $i \neq j$, and complex values $\zeta_k^{(l_k)}$, $k=0,\ldots,m-1$, $l_k=0,\ldots,\nu_k-1$, there exists a unique Laurent polynomial $P_m \in \Lambda_{-p,r-p-1}$, with $$r=\sum_{k=0}^{m-1} \nu_k$$ satisfying
\begin{equation}\label{Hermitecondit}
P_m^{(l_k)}(\tilde{z}_k)=\zeta_k^{(l_k)}, \quad \quad k=0,1,\ldots,m-1, \quad \quad l_k=0,1,\ldots,\nu_k-1.
\end{equation}
\end{theorem}

\begin{DT}
Let us first show uniqueness. If there are two solutions $P_m,Q_m \in \Lambda_{-p,r-p-1}$, then the difference $R_m=P_m-Q_m\in \Lambda_{-p,r-p-1}$ satisfies
$$R_m^{(l_k)}(\tilde{z}_k)=0, \quad \quad k=0,1,\ldots,m-1, \quad \quad l_k=0,1,\ldots,\nu_k-1.$$
Thus, $\tilde{z}_k$ is a zero of $R_m$ of order at least $\nu_k$, and hence $R_m(z)=\frac{S_{r-1}(z)}{z^p}$ where $S_{r-1} \in \PP_{r-1}$ has at least $r$ zeros. It must hold thus $R_m \equiv 0$.

The existence is a consequence of uniqueness. The conditions (\ref{Hermitecondit}) determine a system of $r$ linear equations for the $r$ unknown coefficients $c_j$ of $S_{r-1}(z)=c_0+c_1z+\cdots+c_{r-1}z^{r-1}$, and the matrix of this system is nonsingular, because of the uniqueness of its solution and since the homogeneous problem $\zeta_k^{(l_k)}=0$ admits the solution $R_m \equiv 0$. Hence, (\ref{Hermitecondit}) has a unique solution for arbitrary values $\zeta_k^{(l_k)}$.
\begin{flushright}
$\Box$
\end{flushright}
\end{DT}
An explicit formula of Lagrange type for the Laurent polynomial $P_m$ is given by

\begin{theorem}\label{TheoHermite2}
The unique solution of the problem stated in Theorem \ref{TheoHermite1} is given by
\begin{equation}\label{Hermitesolution}
P_m(z)=\sum_{k=0}^{m-1} \sum_{l_k=0}^{\nu_k-1} \zeta_k^{(l_k)} L_{k,l_k}(z),
\end{equation}
where  $L_{k,l_k}$ are the generalized fundamental Lagrange-Laurent polynomials defined from the auxiliary Laurent polynomials
\begin{equation}\label{Hermitel}
\displaystyle{\ell_{k,l_k}(z)= \frac{(z-\tilde{z}_k)^{l_k}}{l_k!}\cdot \left( \frac{\tilde{z}_k}{z} \right)^p \cdot \prod_{j=0, \; j \neq k}^{m-1} \left( \frac{z-\tilde{z}_j}{\tilde{z}_k - \tilde{z}_j} \right)^{\nu_j} \in \Lambda_{-p,r+l_k-p-\nu_k}, \quad \begin{array}{ccl} k&=&0,1,\ldots,m-1, \\ l_k&=&0,1,\ldots,\nu_k-1, \end{array} }
\end{equation}
as follows.  Setting
\begin{equation}\label{Hermite2}
\displaystyle{L_{k,\nu_k-1}(z)=\ell_{k,\nu_k-1}(z)\in \Lambda_{-p,r-p-1}, \quad \quad k=0,1,\ldots,m-1,}
\end{equation}
let us define recursively for $l_k=\nu_k-2,\nu_k-3,\ldots,0,$
\begin{equation}\label{Hermite3}
\displaystyle{L_{k,l_k}(z)=\ell_{k,l_k}(z) - \sum_{s=l_k+1}^{\nu_k-1} \ell_{k,l_k}^{(s)}(\tilde{z}_k)L_{k,s}(z)\in \Lambda_{-p,r-p-1}.}
\end{equation}
\end{theorem}

\begin{DT}
The result follows from (\ref{Hermitesolution}) if for fixed $k\in \{ 0,1,\ldots,m-1 \}$ and $l_k \in \{ 0,1,\ldots,\nu_k-1 \}$, we prove for all $j=0,1,\ldots,m-1$ and $\sigma=0,1,\ldots,\nu_j-1,$ that
\begin{equation}\label{Hermite4}
L_{k,l_k}^{(\sigma)}(\tilde{z}_j)=\left\{ \begin{array}{cl} 1, \;&\;\textrm{if } j=k \;\textrm{and } l_k=\sigma, \\ 0, \;&\;\textrm{otherwise}.\end{array} \right.
\end{equation}
For fixed $k$, we proceed by backward induction on the index $l_k$. Taking $l_k=\nu_k-1$ in (\ref{Hermitel}), the proof trivially follows from (\ref{Hermite2}) by definition. Suppose that the statement holds for all $l_k \in \{\nu_k-t,\ldots,\nu_k-1\}$, for some $t \in \{2,\ldots,\nu_k-1 \}$ and let us prove (\ref{Hermite4}) for $l_k=\nu_k-t-1$. From (\ref{Hermite3}),
\begin{equation}\label{Hermite5}
\displaystyle{L_{k,\nu_k-t-1}^{(\sigma)}(z)=\ell_{k,\nu_k-t-1}^{(\sigma)}(z) - \sum_{s=\nu_k-t}^{\nu_k-1} \ell_{k,\nu_k-t-1}^{(s)}(\tilde{z}_k)L_{k,s}^{(\sigma)}(z).}
\end{equation}

Take first $z=\tilde{z}_k$. If $\sigma \in \{0,\ldots,\nu_k-t-1\},$ then  from the definition in (\ref{Hermitel}) it follows that $\ell_{k,\nu_k-t-1}^{(\sigma)}(\tilde{z}_k)=\delta_{\sigma,\nu_k-t-1}$ (Kronecker delta symbol) and $L_{k,s}^{(\sigma)}(\tilde{z}_k)=0$ by induction. Thus, $L_{k,\nu_k-t-1}^{(\sigma)}(\tilde{z}_k)=\delta_{\sigma,\nu_k-t-1}$. If $\sigma \in \{\nu_k-t,\ldots,\nu_k-1\},$ then it follows by induction that all the terms under the summation symbol in (\ref{Hermite5}) vanish except when $s=\sigma$. Hence, $$L_{k,\nu_k-t-1}^{(\sigma)}(\tilde{z}_k)=\ell_{k,\nu_k-t-1}^{(\sigma)}(\tilde{z}_k)-\ell_{k,\nu_k-t-1}^{(\sigma)}(\tilde{z}_k)L_{k,\sigma}^{(\sigma)}(\tilde{z}_k)=0.$$

Set finally $z=\tilde{z}_j$, $j \neq k$. From the definition in (\ref{Hermitel}) it follows in this case $\ell_{k,\nu_k-t-1}^{(\sigma)}(\tilde{z}_j)=0$ and by induction, $L_{k,s}^{(\sigma)}(\tilde{z}_j)=0$, for all $\sigma=0,1,\ldots,\nu_j-1$. Hence, $L_{k,\nu_k-t-1}^{(\sigma)}(\tilde{z}_j)=0$ and this completes the proof.
\begin{flushright}
$\Box$
\end{flushright}
\end{DT}

\begin{remark}
In the real line situation only the cases $\nu_k=1$ (see e.g. \cite{Del1,Del2,Del4,Del5,Del8}) and $\nu_k=\nu$ for all $k=0,\ldots,m-1$, see \cite{Del9}, have been considered.
\end{remark}

Notice that we can now estimate the integral (\ref{integralSzego}) by making $I_{\mu}(F) \approx I_{\mu}(L)$. The result is an alternative to the quadrature formula (\ref{Inmu}):
\begin{equation}\label{alternativeqf}
I_{\mu}(F) \approx \sum_{k=0}^{m-1} \sum_{l_k=0}^{\nu_k-1} \omega_{k,l_k} F^{(l_k)}(\xi_k) + \sum_{l=0}^{r-m-1} C_l \eta_l,
\end{equation}
where
\begin{equation}\label{alternativeqf2}
\omega_{k,l_k}=I_{\mu}\left(L_{k,l_k}(z)\right)\quad \quad \textrm{and} \quad \quad \eta_l=I_{\mu}\left( \frac{\omega_m \upsilon_l}{z^p} \right) \quad \quad \textrm{with} \quad \quad \omega_m=\prod_{k=0}^{m-1} (z-\tilde{z}_k),
\end{equation}
for all $k=0,\ldots,m-1$, $l_k=0,\ldots,\nu_k-1$ and $l=0,\ldots,r-m$., $L_{k,l_k}$ being the generalized fundamental Lagrange-Laurent polynomials introduced in Theorem \ref{TheoHermite2}.

We conclude with a result for the particular limit case $r=N$.
\begin{theorem}\label{Caser=N}
In the Lagrange-Laurent interpolation case \eqref{Lagrange-Laurent}, if $r=N,$ then $L$ is the unique Laurent polynomial in $\Lambda_{-p,N-p-1}$ that interpolates $F$ at the set of nodes $Z_N$.
\end{theorem}

\begin{DT}
It is clear from the definition of $L$ in (\ref{L}) that it satisfies the same interpolation conditions as $P_m$ in the set of nodes $\Upsilon_m$. So we need to prove the result for the nodes in $Z_N \backslash \Upsilon_m =\left\{ \hat{z}_s \right\}_{s=1}^{N-m}$. Since $Q_{N-m}\in \PP_{N-m-1}$, whose dimension is equal to the number of the least squares nodes, $Q_{N-m}$ is actually the interpolating polynomial: $Q_{N-m}(\hat{z}_s)=G(\hat{z}_s)$, $s=1,\ldots,N-m$. Thus,
$$\begin{array}{ccl}
L(\hat{z}_s)&=&P_m(\hat{z}_s)+\frac{\omega_m(\hat{z}_s)}{\hat{z}_s^p}Q_{N-m}(\hat{z}_s)=P_m(\hat{z}_s)+\frac{\omega_m(\hat{z}_s)}{\hat{z}_s^p}G(\hat{z}_s) \\ \\
&=&P_m(\hat{z}_s)+\frac{\omega_m(\hat{z}_s)}{\hat{z}_s^p}\left( \frac{\left[F(\hat{z}_s)-P_m(\hat{z}_s)\right]\hat{z}_s^p}{\omega_m(\hat{z}_s)} \right) = F(\hat{z}_s).
\end{array}$$
This completes the proof.
\begin{flushright}
$\Box$
\end{flushright}
\end{DT}

\section{Numerical examples}\label{sec4}

In this section we will consider some illustrative numerical examples about the results obtained throughout the paper for the particular case of the positive measure on the unit circle known as {\em Rogers-Szeg\H{o} weight function}.

Consider the weight function defined on the real line
$$\sigma_\gamma(x)=\sqrt{\frac{\gamma}{\pi}}e^{-\gamma x^2}, \quad \quad \gamma>0, \quad \quad x \in \RR,$$
and the following one associated with $\sigma_\gamma$:
\begin{equation}\label{omega}
\omega(\theta)=\sqrt{\frac{\gamma}{\pi}}\sum_{j=-\infty}^\infty{e^{-\gamma(\theta-2\pi j)^2}}, \quad \quad \theta \in [0,2\pi].
\end{equation}
It is not hard to verify that $\omega$ is efectively a weight function and that
$$\int^{2\pi}_{0}{f(\theta)\omega(\theta)d\theta}=\int^{\infty}_{-\infty}{f(x)\sigma_\gamma(x)dx}$$
for all $2\pi$-periodic function $f$ defined on $\RR$ (see \cite[Theorem 3.1]{RuyRogers}).
The Rogers-Szeg\H{o} weight function is obtained by taking $q=q(\gamma)=e^{-\frac{1}{2\gamma}}$. The resulting normalized weight function (that is, $\mu_0=\int_{0}^{2\pi} \omega(\theta)d\theta=1$)
is given by
\begin{equation}\label{rs}
\omega(\theta,q)=\frac{1}{\sqrt{2\pi\log{(1/q)}}}\sum_{j=-\infty}^\infty{\exp{\left(-\frac{(\theta-2\pi j )^2}{2\log{(1/q)}}\right)}}\,\,,\,\, 0<q<1.
\end{equation}

The motivation to deal with this particular measure for our purposes is twofold. On the one hand, because of its importance in applications, it is consider as the analog on the unit circle of the Gaussian distribution on the real line (see \cite{RuyRogers}, \cite{SimonBk} and references therein). On the other hand, because in the limit cases $q \rightarrow 0^+$ and $q \rightarrow 1^{-},$ the Rogers-Szeg\H{o} weight function tends in the first case to the normalized Lebesgue measure $\frac{d\theta}{2\pi}$, $\theta \in [0,2\pi)$, whereas in the second case to a $\delta$-Dirac distribution at $z=1$. When $q \rightarrow 0^+$, the nodes of the associated para-orthogonal polynomials behave as a uniform distribution of nodes on the unit circle, that is the best scenario in our context since we can take $m\sim N$, and interpolatory-type quadrature formulas on the unit circle will produce good results since they tend to Szeg\H{o} rules. In this case, the construction of the function $L$ in (\ref{L}) will be more focused on an interpolatory process rather than a regression one. However, in the pathological case $q \rightarrow 1^{-}$ the situation is just in the opposite direction. The zeros of para-orthogonal polynomials will  accumulate at $z=1$ and it will be expected to obtain $m \sim 1$. Thus, in the construction of the function $L$ in (\ref{L}), the interpolatory process will not be relevant in this case. The behavior of the parameter $m=m(q)$ along with the numerical estimations of (\ref{integralSzego}) depending on the parameter $q$ is hence an attractive problem to analyze numerically in this section. We point out also that Szeg\H{o} quadrature formulas associated with the Rogers-Szeg\H{o} weight function have been considered in \cite{RuyRogers}. In particular, these have been compared there with interpolatory-type rules based on a uniform distribution of nodes on the unit circle.

Many properties for this weight function are explicitly known (see \cite{RuyRogers} and \cite[Chapter 1.6]{SimonBk}). Let us summarize those that are of  our interest  in the following
\begin{proposition}\label{propertiesRogers}
For the Rogers-Szeg\H{o} weight function \eqref{rs} the following explicit formulas hold.
\begin{enumerate}
\item The sequence $\{\mu_k\}_{k=-\infty}^{\infty}$ of trigonometric moments is given by $$\mu_k=\int_{0}^{2\pi} e^{-ik\theta}\omega(\theta,q)d\theta=q^\frac{k^2}{2}, \quad k\in\ZZ.$$
\item The monic Szeg\H{o} polynomials are given by $\rho_0 \equiv 1$ and for all $n\geq 1$,
\begin{equation}\label{monicrogers}
\rho_n(z)=\sum_{j=0}^{n} (-1)^{n-j} \left[ \substack{n \\ j }
\right]_q q^{\frac{n-j}{2}}z^j,
\end{equation}
where for all $0< q < 1$, the usual $q$-binomial coefficients are defined by
\begin{equation}\label{qbin}
\left[ \substack{n \\ j } \right]_q = \frac{(1-q^n) \cdots (1-q^{n-j+1})}{(1-q) \cdots (1-q^j)}, \quad \quad \textrm{for} \quad 1<j<n \quad \textrm{and} \quad \left[ \substack{n \\ 0 } \right]_q = \left[
\substack{n \\ n } \right]_q  \equiv 1.
\end{equation}
\item The  Verblunsky coefficients \eqref{Verblcoeff} are $\delta_n=(-1)^n q^{n/2}$, for all $n\geq 1$.
\item The invariant para-orthogonal polynomials \eqref{para} are given by
$$B_n(z,\tau_n)=\sum_{j=0}^n{(-1)^{n-j}\left[ \substack{n \\ j } \right]_q q^{\frac{n-j}{2}}\left[1+\tau_n(-1)^nq^{j-\frac{n}{2}}\right]z^j},\quad \tau_n \in \TT.$$
\end{enumerate}
\begin{flushright}
$\square$
\end{flushright}
\end{proposition}

\begin{figure}[H]
\centering
\begin{subfigure}[b]{0.49\linewidth}
\includegraphics[width=\linewidth]{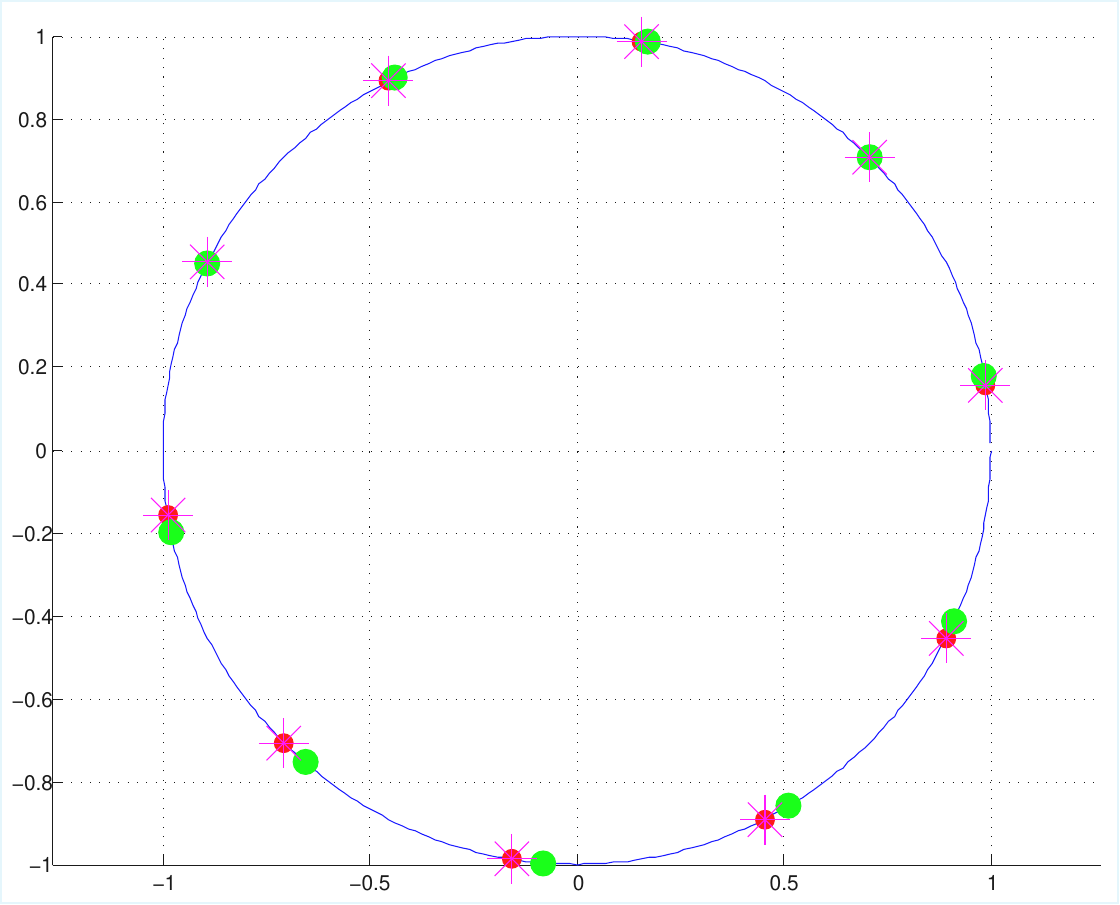}
\caption{Case $q=0.05$ and $N=m=10$.}
\label{fig1a}
\end{subfigure}
\begin{subfigure}[b]{0.49\linewidth}
\includegraphics[width=\linewidth]{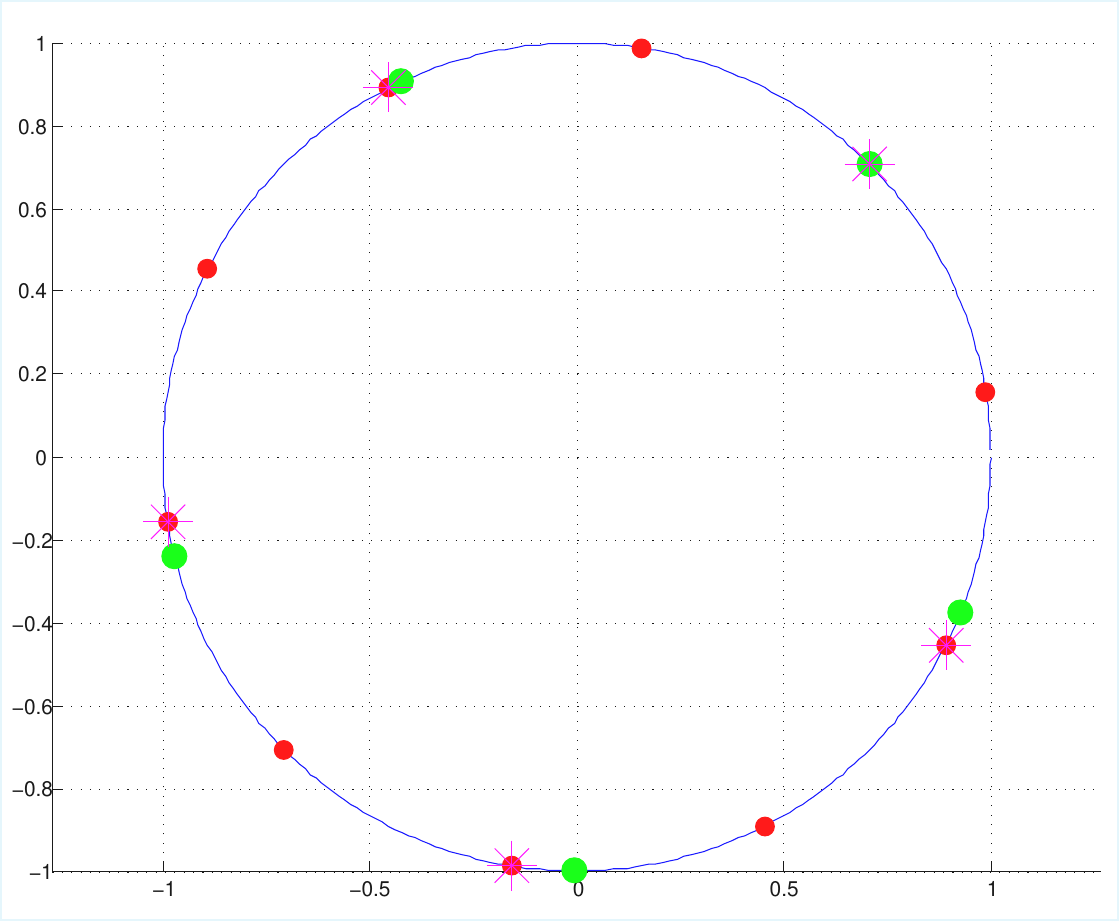}
\caption{Case $q=0.05$, $N=10$ and $m=5$.}
\label{fig1b}
\end{subfigure}
\begin{subfigure}[b]{0.49\linewidth}
\includegraphics[width=\linewidth]{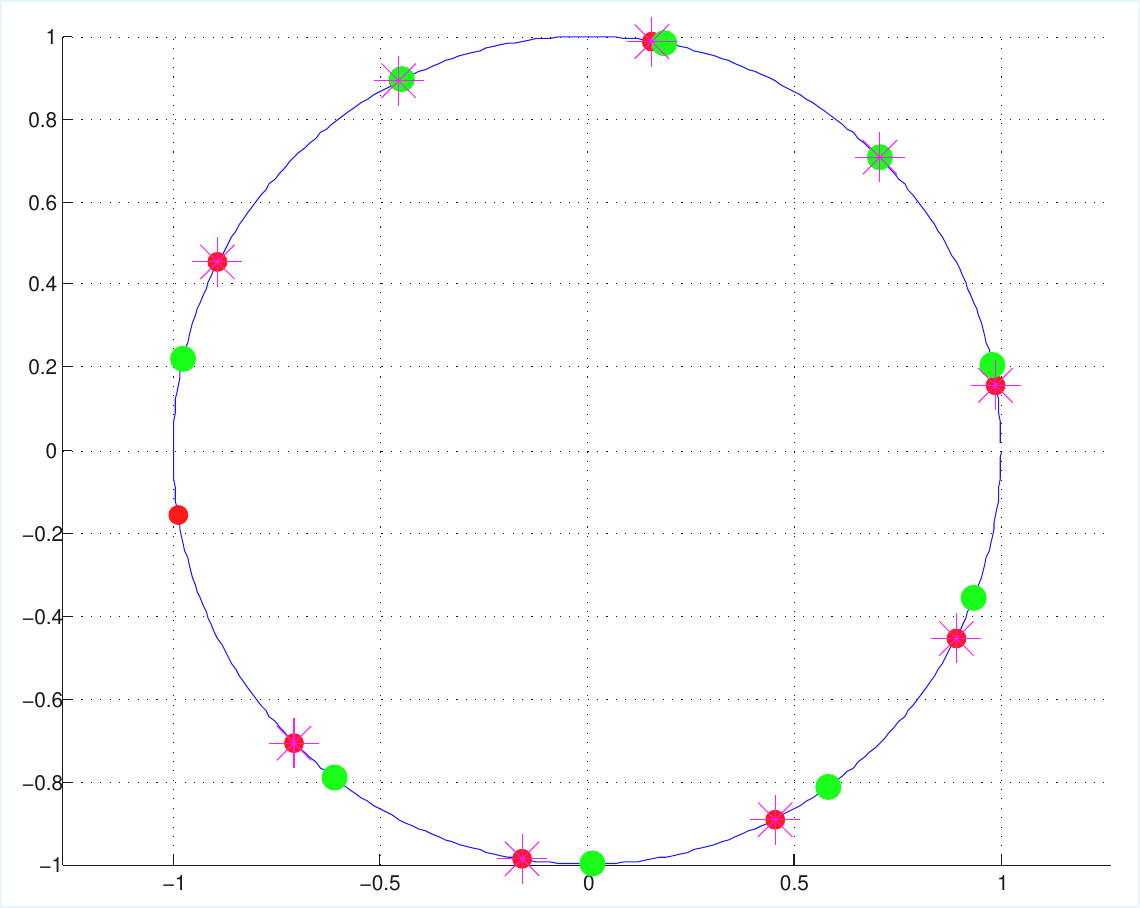}
\caption{Case $q=0.5$, $N=10$ and $m=9$.}
\label{fig1c}
\end{subfigure}
\begin{subfigure}[b]{0.49\linewidth}
\includegraphics[width=\linewidth]{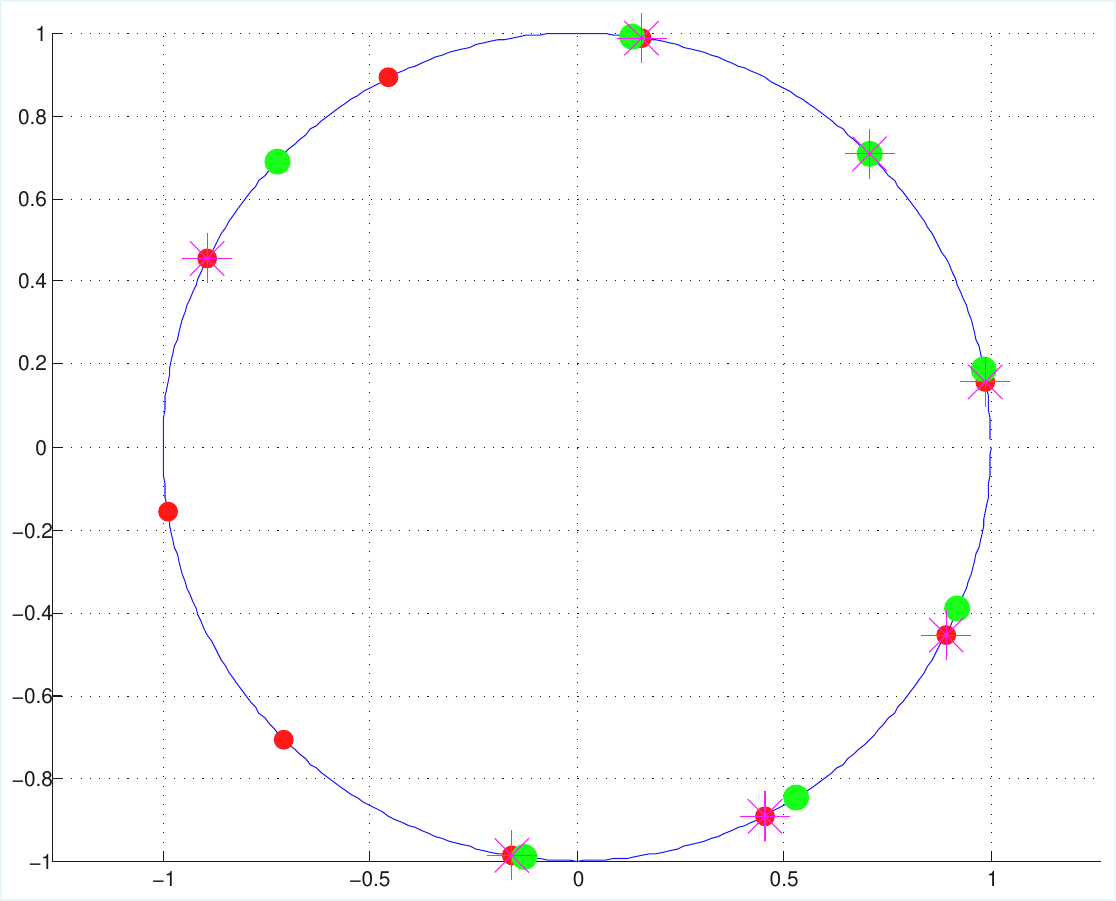}
\caption{Case $q=0.7$, $N=10$ and $m=7$.}
\label{fig1d}
\end{subfigure}
\caption{The set of nodes $Z_{10}$, $\Xi_m$ and $\Upsilon_m$ in some particular cases.}
\label{fig1}
\end{figure}

\begin{figure}[H]
\centering
\begin{subfigure}[b]{0.49\linewidth}
\includegraphics[width=\linewidth]{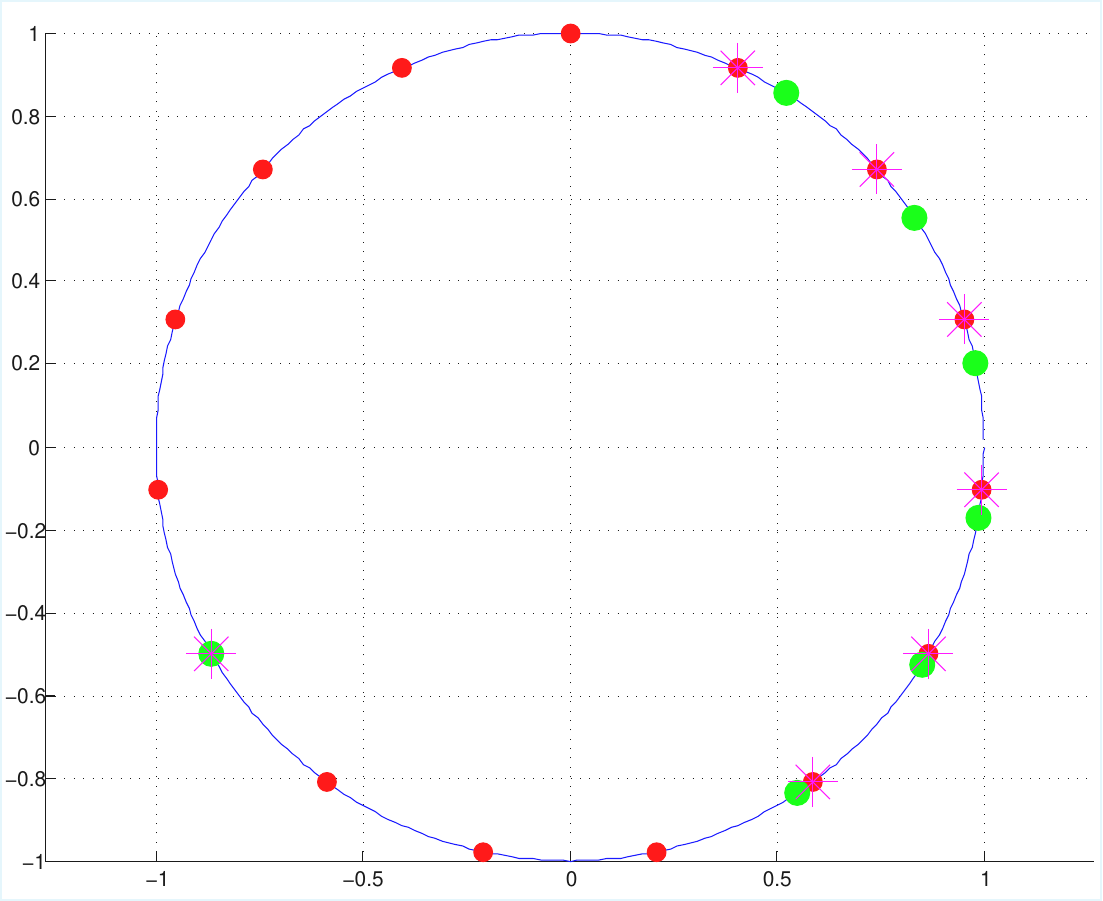}
\caption{Case $q=0.9$, $N=15$, $m=7$ and $\theta_0=\frac{7\pi}{6}$.}
\label{fig2a}
\end{subfigure}
\begin{subfigure}[b]{0.49\linewidth}
\includegraphics[width=\linewidth]{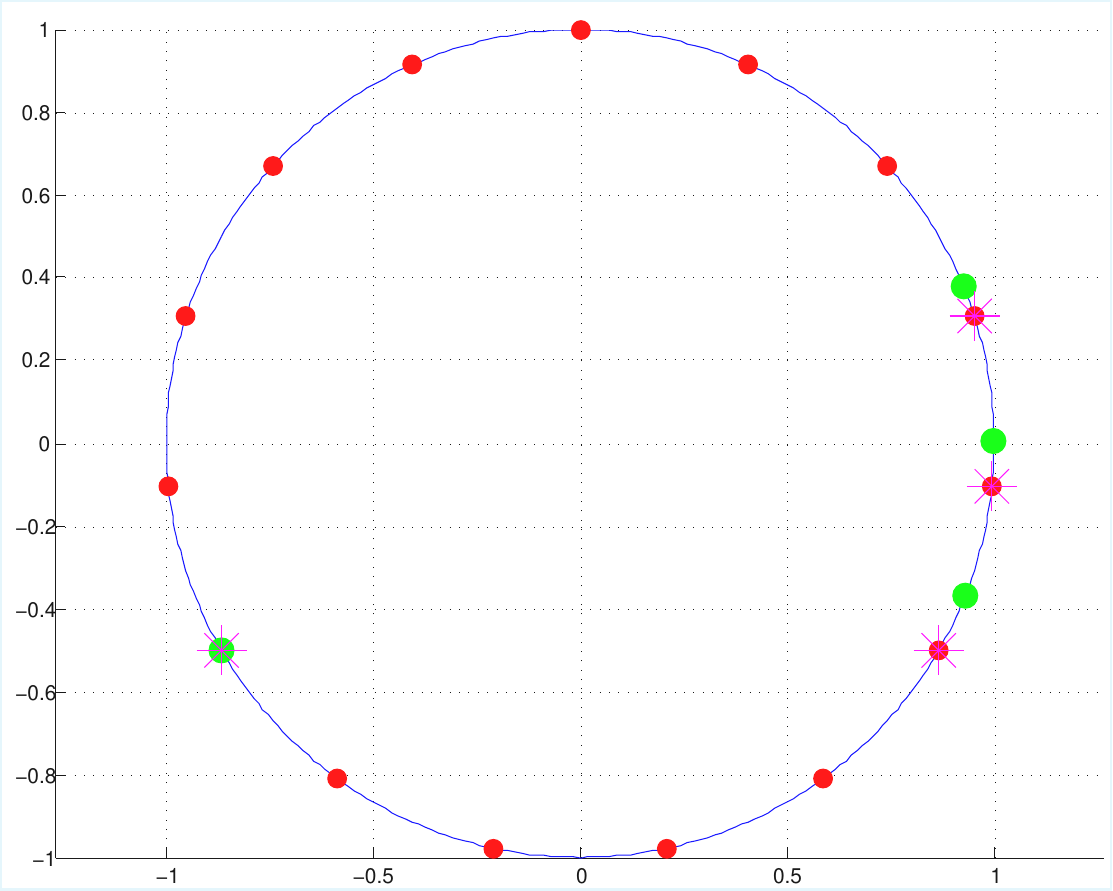}
\caption{Case $q=0.95$, $N=15$, $m=4$ and $\theta_0=\frac{7\pi}{6}$.}
\label{fig2b}
\end{subfigure}
\begin{subfigure}[b]{0.49\linewidth}
\includegraphics[width=\linewidth]{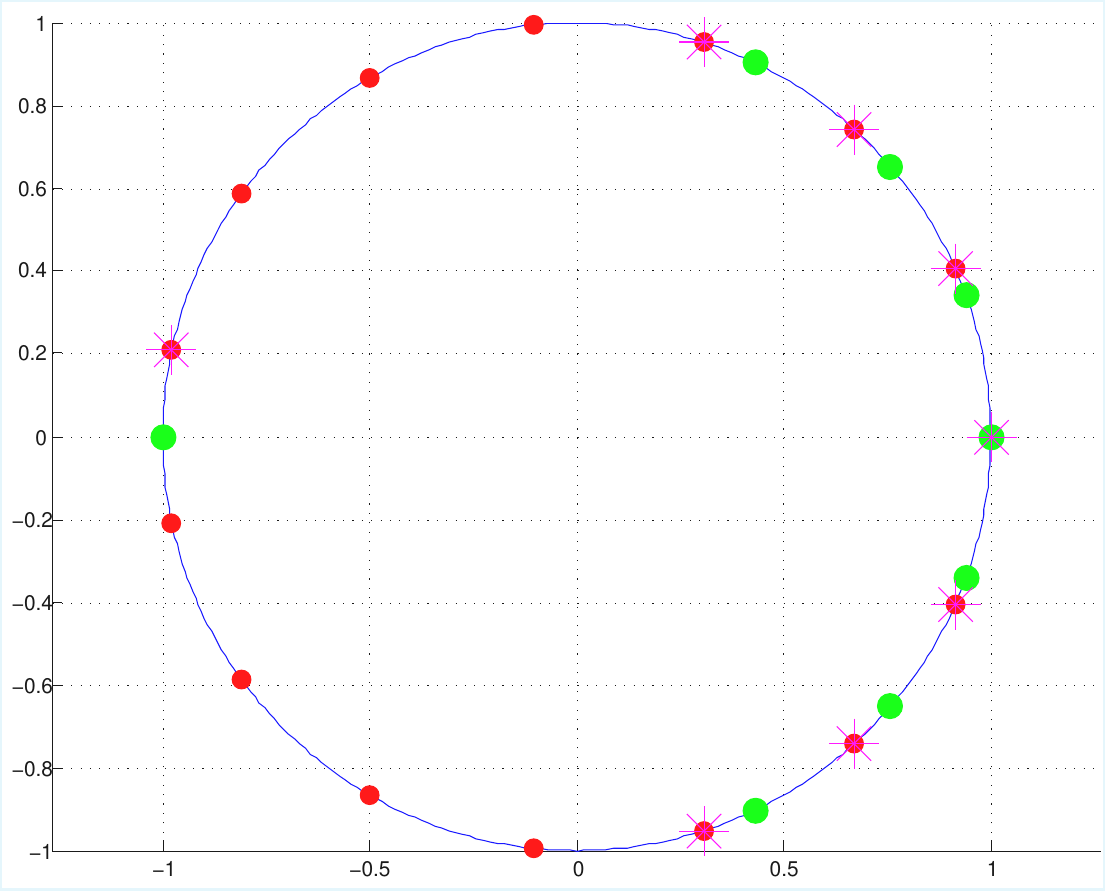}
\caption{Case $q=0.9$, $N=15$, $m=8$ and $\theta_0=0$.}
\label{fig2c}
\end{subfigure}
\begin{subfigure}[b]{0.49\linewidth}
\includegraphics[width=\linewidth]{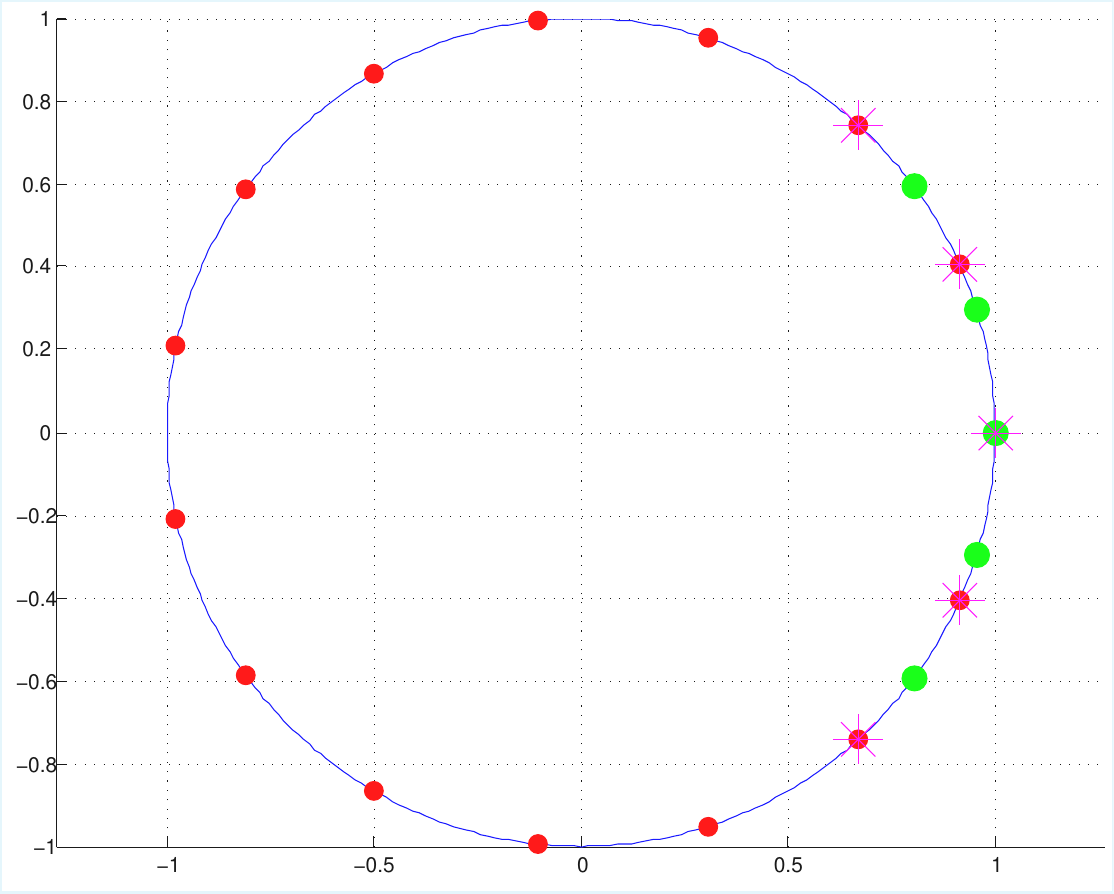}
\caption{Case $q=0.95$, $N=15$, $m=5$ and $\theta_0=0$.}
\label{fig2d}
\end{subfigure}
\caption{The set of nodes $Z_{15}$, $\Xi_m$ and $\Upsilon_m$ in some particular cases for $q \sim 1^{-}$.}
\label{fig2}
\end{figure}

All the computations have been done by making use of \texttt{MATLAB} software. We start by showing in Figures \ref{fig1}-\ref{fig2} the set of nodes $Z_N$, $\Xi_m$ and $\Upsilon_m$ in some particular cases.

The plots in Figure \ref{fig1} show the set of nodes $Z_{10}=\{ z_j \}_{j=0}^{9}$ (red dots, that are uniformly distributed on $\TT$), $\Xi_m=\{ \xi_k \}_{k=0}^{m-1}$ (green dots as the zeros of $B_{m}(z,\tau_m)$) and $\Upsilon_m =\{ \tilde{z}_k \}_{k=0}^{m-1} \subset Z_{10}$ (star-shaped magenta dots, that are the solutions of the problems (\ref{problemk}) for all $k=1,\ldots,m-1$). We have taken different values $1 < m \leq 10$ and $\theta_0=\frac{\pi}{4}$ in all the cases. The parameter $\tau_m$ is computed for each $m$ so that there is a zero of $B_{m}(z,\tau_m)$ at $z_0=\tilde{z}_0=\xi_0=e^{i\frac{\pi}{4}}$ (it is not completely clear in the graphs that this node, of green and star-shaped magenta colours, should also be drawn in red). The Plot (\ref{fig1a}) corresponds to $q=0.05$, $m=10$ and computed value $\tau_{10} \approx 0.283934625753253 - 0.958843641214959i$. This is an extreme situation. Since $q\approx 0$, the nodes of $B_{10}(z,\tau_{10})$ are close to match those of $Z_{10}$ (the Rogers-Szeg\H{o} weight function is close to the Lebesgue measure). Despite we have in mind $m<N$ throughout the paper, in this case $m=N=10$ is the greatest value so that $\Upsilon_{10}$ contains $10$ different nodes. Hence, it has no sense to consider an interpolation-regression model for the estimations of integrals of the form (\ref{integralSzego}) in this case since interpolatory-type rules on $\TT$ are close to the corresponding Szeg\H{o} quadrature formulas (\ref{Inmu}). For the same value $q=0.05$ near zero and $m<10$, the zeros of $B_{m}(z,\tau_m)$ are close to be uniformly distributed on $\TT$ starting at $e^{i\frac{\pi}{4}}$. The distances between two of these consecutive zeros are longer than the fixed distance between two consecutive nodes of $Z_{10}$, and the set $\Upsilon_m$ has thus always $m$ different points. We illustrate this for $m=5$ in the Plot (\ref{fig1b}), where $\tau_5 = 0.477418458445432 + 0.878676058360296i$. Observe in this case that the zeros of $B_{5}(z,\tau_{5})$ are close to each two consecutive points nodes of $Z_{10}$. When the value of $q$ increases, the zeros of the para-orthogonal polynomials tend to move away from the point $z=-1$ since the weight function tends for $q\rightarrow 1^{-}$ to a $\delta$-Dirac distribution at $z=1$. Two results for intermediate values for $q$ and the maximum value that $m$ can take in each case are shown in Plots (\ref{fig1c})-(\ref{fig1d}). The corresponding parameters are $\tau_9=0.665584145627702 - 0.746322815602633i$ and $\tau_7=0.604589851170908 - 0.796536949463826i$, respectively.

The Plots \ref{fig2a}-\ref{fig2b} in Figure \ref{fig2} show the same sets of nodes as in Figure \ref{fig1}, now taking $N=15$, values $q \rightarrow 1^{-}$, the maximum values for the parameter $m$ in each case and $\theta_0=\frac{7\pi}{6}$, which is unnatural since the zeros of $B_{m}(z,\tau_m)$ tend to move away from $z=-1$. It is clearly observed that the value $\frac{N}{m}\geq 1$ increases as $q \rightarrow 1^{-}$. The corresponding $\tau_m$ parameters are $\tau_7=0.993180485559280 + 0.116586976563560i$, $\tau_4=-0.999590226967595 - 0.028624782110457i$, $\tau_8=\tau_5=-1$, respectively. As expected, the freedom of the parameter $\theta_0$ can be used to increase the value of $m$ if some information about the location of the zeros of the para-orthogonal polynomials is known in advance. In the cases (\ref{fig2c})-(\ref{fig2d}) the parameter $\tau_m$ is real because we have forced the zero that we want to leave fixed to be real also (here $\theta_0=0$, the most natural choice). Since the Szeg\H{o} polynomials have real coefficients (notice that $\omega(\theta,q)=\omega(-\theta,q)$ in (\ref{rs})), this implies that the para-orthogonal polynomials will have in those cases also real coefficients, and their zeros appear in complex conjugate pairs on the unit circle. This justifies that depending on whether $m$ is even or odd, the node $-1$ may be part of the sets $\Xi_m$.

In Table \ref{Tab1} we show some relations between the parameters $m$ and $q$ for $N$ fixed. As it was indicated in Section \ref{sec2}, it seems that it is not possible to determine a general formula for $m=m(N)$ for arbitrary measures supported on $\TT$. Based on the results, it is striking that for $q\leq 0.4$ the optimal choice is $m=N$ in all the cases. This is because for these values of $q$, the zeros of the corresponding para-orthogonal polynomials are still close to the uniform distribution that is obtained by making $q\rightarrow 0^{+}$. On the other hand, the situation changes radically when $q$ approaches the limit value $1^{-}$, showing that the value of $m$ decreases rapidly in those cases. If we compare the results in column $N=15$ with those in Figure \ref{fig2}, then  we observe that the value of $m$ can increase a bit when $\theta_0=0$ is chosen, probably due to the symmetric distribution of the zeros of para-orthogonal polynomials.

\begin{table}[H]
\centering
\begin{tabular}{cc|c|c|c|c|c|c|c|c|c|l}
\cline{3-11}
& & \multicolumn{9}{ c| }{Values of $N$} \\ \cline{3-11}
& & $\;$ {\bf 5} $\;$ & $\;${\bf 10}$\;$ & $\;${\bf 15}$\;$ & $\;${\bf 20}$\;$ & $\;${\bf 25}$\;$ & $\;${\bf 30}$\;$ & $\;${\bf 40}$\;$ & $\;${\bf 50}$\;$ & $\;${\bf 100}$\;$ \\ \cline{1-11}
\multicolumn{1}{ |c  }{} &
\multicolumn{1}{ |c| }{$\;${\bf 0.1}$\;$} & 5 & 10 & 15 & 20 & 25 & 30 & 40 & 50 & 100 & \\ \cline{2-11}
\multicolumn{1}{ |c  }{} &
\multicolumn{1}{ |c| }{$\;${\bf 0.2}$\;$} & 5 & 10 & 15 & 20 & 25 & 30 & 40 & 50 & 100 &  \\ \cline{2-11}
\multicolumn{1}{ |c  }{} &
\multicolumn{1}{ |c| }{$\;${\bf 0.3}$\;$} & 5 & 10 & 15 & 20 & 25 & 30 & 40 & 50 & 100 &  \\ \cline{2-11}
\multicolumn{1}{ |c  }{} &
\multicolumn{1}{ |c| }{$\;${\bf 0.4}$\;$} & 5 & 10 & 15 & 20 & 25 & 30 & 40 & 50 & 100 &  \\ \cline{2-11}
\multicolumn{1}{ |c  }{\multirow{2}{*}{Values of $q$}} &
\multicolumn{1}{ |c| }{$\;${\bf 0.5}$\;$} & 4 & 9 & 14 & 19 & 24 & 29 & 39 & 49 & 99 &  \\ \cline{2-11}
\multicolumn{1}{ |c  }{} &
\multicolumn{1}{ |c| }{$\;${\bf 0.6}$\;$} & 4 & 9 & 14 & 19 & 24 & 29 & 39 & 49 & 99 &  \\ \cline{2-11}
\multicolumn{1}{ |c  }{} &
\multicolumn{1}{ |c| }{$\;${\bf 0.7}$\;$} & 3 & 8 & 13 & 18 & 23 & 28 & 38 & 48 & 98 &  \\ \cline{2-11}
\multicolumn{1}{ |c  }{} &
\multicolumn{1}{ |c| }{$\;${\bf 0.8}$\;$} & 3 & 6 & 11 & 15 & 21 & 26 & 36 & 46 & 96 &  \\ \cline{2-11}
\multicolumn{1}{ |c  }{} &
\multicolumn{1}{ |c| }{$\;${\bf 0.9}$\;$} & 2 & 4 & 7 & 11 & 15 & 20 & 29 & 40 & 90 &  \\ \cline{2-11}
\multicolumn{1}{ |c  }{} &
\multicolumn{1}{ |c| }{$\;${\bf 0.95}$\;$} & 1 & 2 & 3 & 6 & 9 & 12 & 20  & 30 & 59 &  \\ \cline{1-11}
\end{tabular}
  \caption{The relation $m=m(N,q)$ for some particular values of $N$ and $q$, and fixed value $\theta_0=\frac{\pi}{6}$.}
  \label{Tab1}
\end{table}

\begin{figure}[H]
\centering
\begin{subfigure}[b]{0.49\linewidth}
\includegraphics[width=\linewidth]{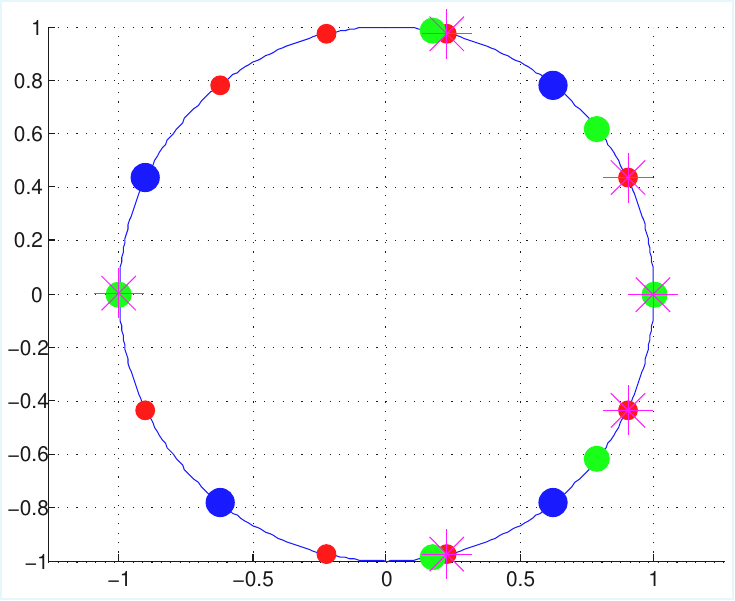}
\caption{Case $q=0.7$, $N=14$, $m=6$ and $\theta_0=0$. Here, $\tau_6=-1$ and $K_{\tilde{\Delta}_4}=\frac{4}{3}$.}
\label{fig3a}
\end{subfigure}
\begin{subfigure}[b]{0.49\linewidth}
\includegraphics[width=\linewidth]{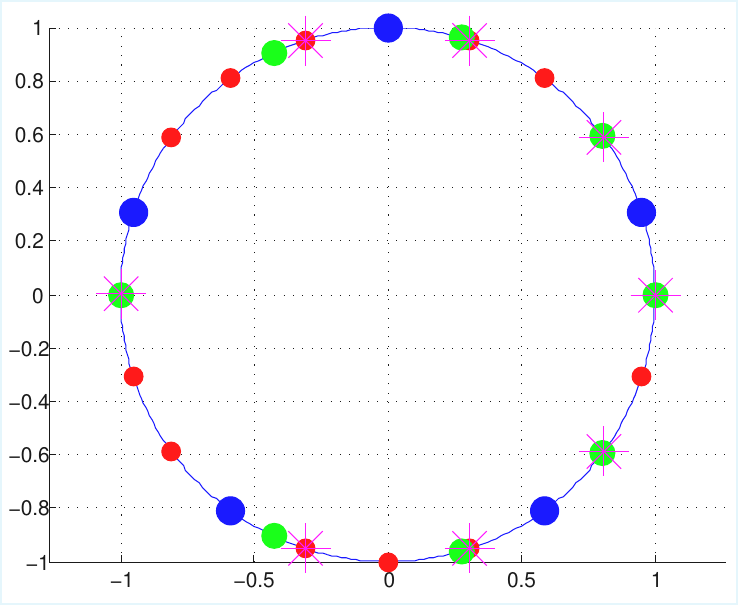}
\caption{Case $q=0.5$, $N=20$, $m=8$ and $\theta_0=0$. Here, $\tau_8=-1$ and $K_{\tilde{\Delta}_5}=1$.}
\label{fig3b}
\end{subfigure}
\begin{subfigure}[b]{0.49\linewidth}
\includegraphics[width=\linewidth]{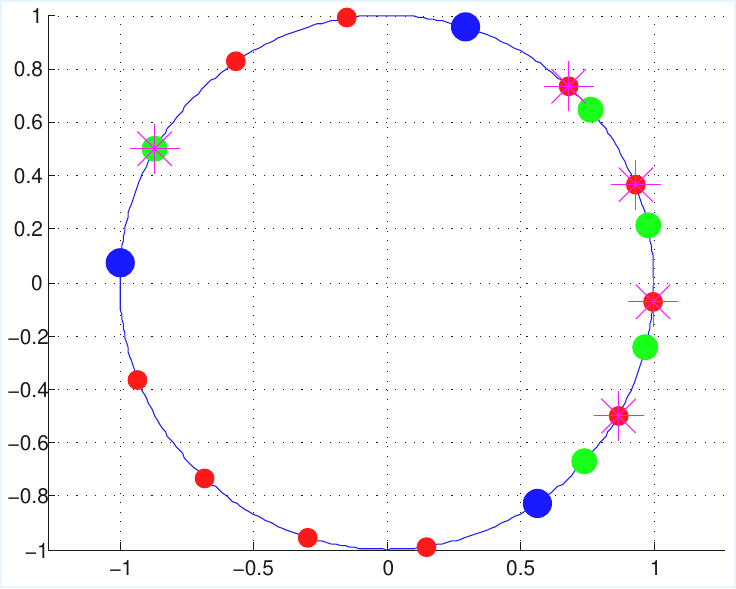}
\caption{Case $q=0.9$, $N=14$, $m=5$ and $\theta_0=\frac{5\pi}{6}$. Here, $\tau_5 \approx 0.996896189168107 - 0.078727301631046i$ and $K_{\tilde{\Delta}_3}=\frac{5}{4}$.}
\label{fig3c}
\end{subfigure}
\begin{subfigure}[b]{0.49\linewidth}
\includegraphics[width=\linewidth]{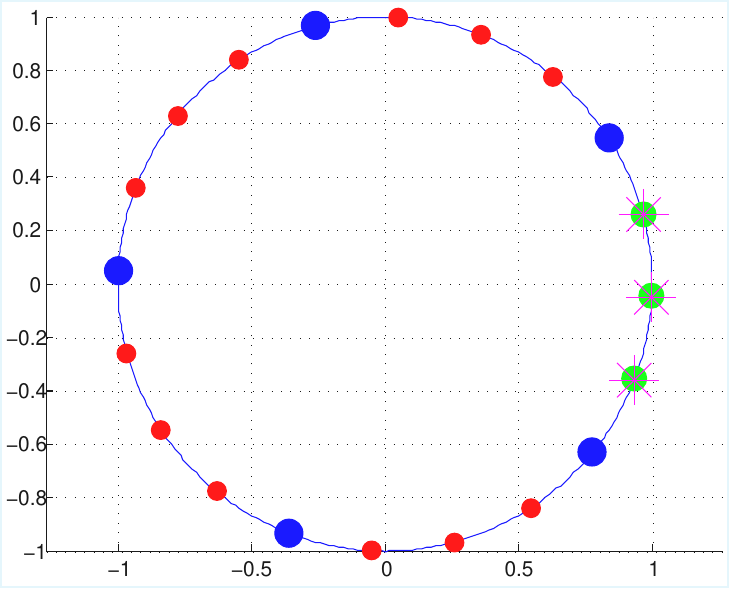}
\caption{Case $q=0.97$, $N=20$, $m=3$ and $\theta_0=\frac{\pi}{12}$. Here, $\tau_3 \approx 0.826882488545677 + 0.562374741730556i$ and $K_{\tilde{\Delta}_5}=1$.}
\label{fig3d}
\end{subfigure}
\caption{The optimal partition $\tilde{\Delta}_k$ for some particular choices of the parameters $m$,$q$ and $N$, along with the computed parameter $K_{\tilde{\Delta}_k}$ in (\ref{Kopt}) that measures the deviation of $\tilde{\Delta}_k$ from the uniform distribution. }
\label{fig3}.
\end{figure}

To the points represented in the Figures \ref{fig1}-\ref{fig2}, we include in Figure \ref{fig3} in blue dots the set of points of the discarded nodes (red dots) that constitute the subpartition $\tilde{\Delta}_k$, of all possible subpartitions formed by those points, which is closest to being uniform. The parameter $K_{\tilde{\Delta}_k}$ in (\ref{Kopt}) that measures the deviation of $\tilde{\Delta}_k$ from the uniform distribution has been computed in each case. In Plot \ref{fig3a} the subpartition $\tilde{\Delta}_4$ contains four points that are placed on the unit circle forming a polygon that is close to being a square. No other combination of these $8$ discarded points will provide a polygon with more vertices on the unit circle that is closer to being regular. In this case, another possible solution with four vertices was found. The optimal situation is achieved in Plot \ref{fig3b}, a very exceptional situation. In this case, the solution is unique. Since $\theta_0=0$ and $\tau_8\in \RR$, the para-orthogonal polynomial $B_8(z,-1)$ has real coefficients and, therefore, its zeros are real ($B_8(\pm 1,-1)=0$) or appear in complex conjugate pairs on $\TT$. The other $10$ discarded nodes on $\TT \backslash \{ \pm 1 \}$ are placed so that a regular pentagon can be formed from $5$ of them. As indicated before, the developed theory is still valid if we initially choose an arbitrary sequence of nodes $Z_N$ not necessarily equally spaced. The advantage of considering a uniformly distributed sequence of points is that, in this case, it is reasonably easy to obtain regular polygons inscribed in the unit circle using the discarded nodes (that is, the optimal value $K_{\tilde{\Delta}_k}=1$ is achieved). However, these regular polygons have a small number of vertices, in many cases, no more than five. Two more examples have been included in Plots \ref{fig3c}-\ref{fig3d}, with parameter $q\rightarrow 1^{-}$, so the zeros of para-orthogonal polynomials tend to accumulate at $z=1$. In Plot \ref{fig3c} it is possible to obtain from the discarded nodes (which except for the zero at $e^{\frac{5\pi}{6}i}$ that has been fixed, are all at the left of the zeros of $B_5(z,\tau_5)$) a triangle that is close to being equilateral. In Plot \ref{fig3d}, since $m$ takes a small value, it is possible to obtain a regular pentagon from the discarded nodes.

Next some examples involving the numerical approximation of integrals are presented. We start by computing the following three quadrature formulas:
\begin{enumerate}
\item An interpolatory-type rule based on the $N$ equally spaced nodes $Z_N$. As the validity domain we chose linear subspaces of Laurent polynomials of the form $\Lambda_{-\ell,\ell}$ ($N=2\ell+1$) or $\Lambda_{-\ell,\ell-1}$ ($N=2\ell$) i.e., the CMV ordering. In this case, an explicit expression for the weights can be found in \cite[Corollary 3.6]{RuyRogers}: if the $N$-point interpolatory-type rule
\begin{equation}\label{INRogers}
I_N^{\omega}(F)=\sum_{j=0}^{N-1} \lambda_j F(z_j)
\end{equation} is exact in $\Lambda_{-r,s}$ with $N=r+s+1$ and the set of nodes $Z_N$ are the zeros of $z^N-\tau$ with $\tau \in \TT$, then
\begin{equation}\label{pesosRogersinterpolat}
\lambda_j=\frac{z_j^r}{N\tau}\sum_{k=1}^{N} q^{\frac{(k-s-1)^2}{2}}z_j^k, \quad \quad j=0,\ldots,N-1.
\end{equation}
Recall that these weights are in general complex numbers, as shown in Table \ref{Tab2} for particular values of $\theta_0$, $q$ and $N$.
\item An interpolatory-type rule
\begin{equation}\label{ImRogers}
I_m^{\omega}(F)=\sum_{k=0}^{m-1} \tilde{\lambda}_k F(\tilde{z}_k)
\end{equation}
based on the $m$ distinct nodes $\Upsilon_m$. The domain of validity is now of the form $\Lambda_{-\tilde{\ell},\tilde{\ell}}$ ($m=2\tilde{\ell}+1$) or $\Lambda_{-\tilde{\ell},\tilde{\ell}-1}$ ($m=2\tilde{\ell}$). The weights can not be computed in this case from (\ref{pesosRogersinterpolat}) since $\Upsilon_m$ is not a set of uniformly distributed nodes. The quadrature formula can be obtained thus from the direct integration of the Laurent polynomial of Lagrange interpolation $P_m$ given by (\ref{Lagrange-Laurent}). However, the computational effort is extremely reduced if we take into account that an explicit expression for the trigonometric moments is known (Theorem \ref{propertiesRogers}-1). Recall that the parameter $m$ is chosen so that it is the largest possible for which set $\Upsilon_m$ can be constructed. Some illustrative examples are presented in Table \ref{Tab3}. Here we must emphasize that these $m$ points mimic the zeros of $B_m(z,\tau_m)$, so these interpolatory-type rules are close to be Szeg\H{o} rules, that are exact in the largest linear subspace of Laurent polynomials (\ref{SantosNjastad}) with $n$ replaced by $m$. Thus, for example, it follows from Table \ref{Tab1} that if $N=25$ and $q=0.8$ we have to take $m=21$. The interpolatory-type quadrature formulas based on the nodes $Z_{25}$ and $\Upsilon_{21}$ are exact in $\Lambda_{-12,12}$ and $\Lambda_{-10,10}$, respectively. But the second one is close to be a $21$-point Szeg\H{o} rule, that is exact in a linear subspace of Laurent polynomials that contains $\Lambda_{-20,20}$. This means that we expect good results for this interpolatory-type quadrature formula, despite the fact we are using the information of the integrand provided by only a subset of the equally spaced nodes.
\item The new numerical procedure (\ref{alternativeqf})-(\ref{alternativeqf2}), first for the Lagrange interpolation case. The rule is obtained from the direct integration of the Laurent polynomial $L$ given by (\ref{L}), computed as described in Section \ref{sec2}. Here we can make use again of the explicit expressions for the trigonometric moments to reduce the computational effort. Since $L$ is constructed from interpolation on the nodes $\Upsilon_m$, it has at least the same domain of validity as the interpolatory-type rule based on these $m$ points described in the previous item. The rest of the information on the nodes $Z_N \backslash \Upsilon_m$ has been used to improve the accuracy of the approximation by a process of simultaneous complex regression, taking into account that the resulting $L$ belongs to a linear subspace of Laurent polynomials of dimension $r$, $m \leq r \leq N$. So, we expect better results than the two previous rules.
\end{enumerate}

We present some first results in Table \ref{Tab2}. Unlike Szeg\H{o} quadrature formulas, interpolatory-type rules do not have positive weights; in general  they are complex numbers. For that reason, the convergence of the process is not always assured for integrands that are Riemann-Stieltjes integrable with respect to $\omega(\theta,q)$. We compute these values in some particular cases. The weights correspond to the nodes ordered counterclockwise starting from $e^{i\theta_0}$. Since the rules are exact for the constant function $F\equiv 1$, the sum of the weights is always equal to $\mu_0=1$. A sufficient condition to ensure the real character of the weights can be found in \cite[Corollary 3.8]{RuyRogers} (for measures having real trigonometric moments). But even so, negative weights may appear, and this is the case of Column 2 (something that becomes more pronunciated as the parameter $q$ approaches the value $1$). These sufficient conditions are satisfied here because we have taken $\theta_0=0$ (so, the nodes are real or appear in complex conjugate pairs) and the interpolatory-type rule is exact in a linear subspace of Laurent polynomials of the form $\Lambda_{-r,s}$ with $|r-s| \leq 1$. The conditions of that result are clearly not necessary, se e.g. Column 4.

\begin{table}[H]
\centering
\resizebox{1.00\textwidth}{!}{
\begin{tabular}{|c|c|c|c|}
\hline
$\theta_0=\frac{\pi}{3}$, $q=0.5$, $N=10$ & $\theta_0=0$, $q=0.8$, $N=12$ & $\theta_0=\frac{4\pi}{3}$, $q=0.95$, $N=8$ & $\theta_0=\frac{4\pi}{3}$, $q=0.95$, $N=9$ \\ \hline
$\begin{array}{c}
0.136489850028398 - 0.000014950498919i \\
0.039746488789610 + 0.000014950498919i \\
0.006539440382177 - 0.000014950498919i \\
0.000709671662507 + 0.000014950498919i \\
0.001457311373535 - 0.000014950498919i \\
0.012728899971602 + 0.000014950498919i \\
0.063866165914576 - 0.000014950498919i \\
0.181479357664911 + 0.000014950498919i \\
0.291690390674189 - 0.000014950498919i \\
0.265292423538495 + 0.000014950498919i
\end{array}$ &
$\begin{array}{c}
  0.439838082705365 \\
  0.241412252032503 \\
  0.036118270265519 \\
  0.003127493457543 \\
 -0.001077897180606 \\
  0.000956654882583 \\
 -0.000911629620451 \\
  0.000956654882583 \\
 -0.001077897180606 \\
  0.003127493457543 \\
  0.036118270265519 \\
  0.241412252032503 \\
\end{array}$
&
$\begin{array}{c}
0.047359416517110 - 0.071817368360745i \\
-0.106206727880972 + 0.071817368360745i \\
0.407373010746110 - 0.071817368360745i \\
0.737567691449938 + 0.071817368360745i \\
-0.105911970428243 - 0.071817368360745i \\
0.048340318334377 + 0.071817368360745i \\
-0.014675564657243 - 0.071817368360745i \\
-0.013846174081077 + 0.071817368360745i \\
\end{array}$ &
$\begin{array}{c}
0.005240568499150 \\
-0.015011776306052 \\
0.085113936805406 \\
0.852107190422470 \\
0.085113936805407 \\
-0.015011776306052 \\
0.005240568499151 \\
-0.001396324209739 \\
-0.001396324209739 \\
\end{array}$
\\ \hline
\end{tabular}
}
\caption{Weighs of the interpolatory-type quadrature formulas $I_N^{\omega}(F)$ given by (\ref{INRogers}) in some particular cases (equally spaced nodes).}
\label{Tab2}
\end{table}

Table \ref{Tab3} shows the angles of the nodes ordered counterchockwise starting from $\theta_0$. In the first column, the nodes $Z_N$ marked with $(\checkmark)$ are those that belong to the set $\Upsilon_m$. In the second column, the zeros of $B_m(z,\tau_m)$ are shown in the same row as their corresponding value from set $\Upsilon_m$ i.e., the solutions of the problems (\ref{problemk}) (in Case 1,  $\tau_9\approx 0.961720991826004 - 0.274030534578186i$, and in Case 2, $\tau_4\approx-0.998305536515183 - 0.058189825254356i$). In the third column we have computed the corresponding quadrature weights. In Case 1, all of them are positive. The reason for this might be that this interpolatory-rule is close to be a Szeg\H{o} quadrature formula, so it is reasonable that the weights approach positive values. But this is not true in general, as shown in Case 2.

\begin{table}[H]
\centering
\resizebox{1.00\textwidth}{!}{
\begin{tabular}{c|c|c|c|}
\cline{2-4}
 & Nodes $Z_N$ and $\Upsilon_m$ & Zeros of $B_m(z,\tau_m)$ & Weights \\ \cline{1-4}
\multicolumn{1}{ |c| }{$\begin{array}{c} \textbf{Case 1}: \\ \\ \begin{array}{l} N=15, \\ m=9, \\ q=0.85, \\ \theta_0=\frac{5\pi}{6}. \end{array} \end{array}$} & $\begin{array}{l}\theta_0\approx2.617993877991494 \;(\checkmark) \\ $\quad \quad$3.036872898470134 \\ $\quad \quad$3.455751918948772 \\ $\quad \quad$3.874630939427412 \\ $\quad \quad$4.293509959906051 \\ $\quad \quad$4.712388980384690 \;(\checkmark) \\ $\quad \quad$5.131268000863329 \;(\checkmark) \\ $\quad \quad$5.550147021341968 \;(\checkmark) \\ $\quad \quad$5.969026041820607 \;(\checkmark) \\ $\quad \quad$0.104719755119660 \;(\checkmark) \\ $\quad \quad$0.523598775598298 \;(\checkmark) \\ $\quad \quad$0.942477796076938 \;(\checkmark) \\ $\quad \quad$1.361356816555578 \;(\checkmark) \\ $\quad \quad$1.780235837034217 \\ $\quad \quad$2.199114857512855 \end{array}$  & $\begin{array}{l}\theta_0 \approx 2.617993877991494  \\ \\ \\ \\ \\ $\quad \quad$4.791809174919338 \\ $\quad \quad$5.264512395959146 \\ $\quad \quad$5.676656940456297 \\ $\quad \quad$6.065035322017836 \\ $\quad \quad$0.163174054491049 \\ $\quad \quad$0.550360774964562 \\ $\quad \quad$0.959655336539369 \\ $\quad \quad$1.426042704164285 \\ \\ \\ \end{array}$ & $\begin{array}{l} 0.400459456060581  \\ \\ \\ \\ \\ 0.306250761095405 \\ 0.178586805292540 \\ 0.079176983347153\\ 0.000198428331037 \\ 0.000003089382994  \\ 0.007062147533195 \\ 0.001462718217120 \\ 0.026799610739976 \\ \\ \\ \end{array}$ \\ \cline{1-4}
\multicolumn{1}{ |c| }{$\begin{array}{c} \textbf{Case 2}: \\ \\ \begin{array}{l} N=9, \\ m=4, \\ q=0.9, \\ \theta_0=-\frac{\pi}{6}. \end{array} \end{array}$} & $\begin{array}{l}\theta_0\approx5.759586531581287 \;(\checkmark) \\ $\quad \quad$0.174532925199433\;(\checkmark)  \\ $\quad \quad$0.872664625997165\;(\checkmark) \\ $\quad \quad$1.570796326794897 \\ $\quad \quad$2.268928027592629 \\ $\quad \quad$2.967059728390360 \\ $\quad \quad$3.665191429188092\;(\checkmark) \\ $\quad \quad$4.363323129985823 \\ $\quad \quad$5.061454830783555 \end{array}$ & $\begin{array}{l}\theta_0 \approx 5.759586531581287 \\ $\quad \quad$0.013932415392792 \\ $\quad \quad$0.552145230405753 \\ \\ \\ \\ $\quad \quad$3.640310346572479 \\ \\ \\ \end{array}$ & $\begin{array}{r}0.018604143064766 - 0.013341862967978i \\ -0.001616984685047 + 0.003439169388233i \\ 0.275749994559734 - 0.015171806777346i  \\ \\ \\ \\ 0.707262847060548 + 0.025074500357091i \\ \\ \\ \end{array}$ \\ \cline{1-4}
\end{tabular}}
\caption{Nodes and weights of the interpolatory-type quadrature formula $I_m^{\omega}(F)$ given by (\ref{ImRogers}).}
\label{Tab3}
\end{table}

Some interesting consequences can be derived from the results presented in Table \ref{Tab2}. We see that in the case of the CMV ordering, the weights $\lambda_j$ in the rule
 $I_N^{\omega}(F)$ given by (\ref{INRogers}) are real when $N$ is odd and if $N$ is even, the weights are complex quantities, but all of them have equal imaginary part, up to a sign, independently of the value $\theta_0 \not\in \{0,\pi \}$. For $\theta_0 \in \{0,\pi \}$, the weights are also real if $N$ is even. This property is particular for the Rogers-Szeg\H{o} case, and even for this weight function, it is not verified if we consider interpolatory-type quadrature formulas based on uniformly distributed nodes on the unit circle but exact in linear subspaces of the form $\Lambda_{-r,s}$ with $|r-s| > 1$, see e.g. \cite[Table 3]{RuyRogers}. We can prove this result, that is new in the literature, and that has repercussions on the computational effort required to compute these formulas.
\begin{proposition}
Let $I_N^{\omega}(F)$ be the $N$-point interpolatory-type quadrature formula for the Rogers-Szeg\H{o} weight function $\omega$, based on the nodes $Z_N$ for arbitrary $\theta_0$ (see \eqref{INRogers}, \eqref{rs} and \eqref{ZN}, respectively). Suppose that $I_N^{\omega}(F)$ is exact in $\Lambda_{-\ell,\ell},$ if $N=2\ell+1,$ and either $\Lambda_{-\ell,\ell-1}$ or $\Lambda_{-(\ell-1),\ell},$ if $N=2\ell$. Then, the weights $\lambda_j$ are
\begin{equation}\label{pesosRogersinterpolatexplicit}
\lambda_j=\frac{1}{N} \times \left\{ \begin{array}{ccl} 1+ 2\cdot\displaystyle{\sum_{k=1}^{\ell} } q^{\frac{k^2}{2}}\cos \left[\left( \theta_0+hj \right)k\right], &\textrm{if}& N=2\ell+1, \\ \\ 1+ (-1)^j q^{\ell^2/2}\left[ \cos(\ell\theta_0)+i\epsilon\sin(\ell\theta_0) \right] + 2\cdot\displaystyle{\sum_{k=1}^{\ell-1} } q^{\frac{k^2}{2}}\cos \left[\left( \theta_0+hj \right)k\right], &\textrm{if}& N=2\ell,
\end{array}\right.
\end{equation}
where $$\epsilon=\left\{\begin{array}{rcl} 1, &\textrm{if}& I_{2\ell}^{\omega}(F) \;\textrm{is exact in } \Lambda_{-\ell,\ell-1}, \\ \\ -1, &\textrm{if}& I_{2\ell}^{\omega}(F) \;\textrm{is exact in } \Lambda_{-(\ell-1),\ell} . \end{array}\right.$$
\end{proposition}

\begin{DT}
Suppose first the case $N=2\ell+1$, so $r=s=\ell$. Rearranging the sum in (\ref{pesosRogersinterpolat}) where $\tau=e^{i(2\ell+1)\theta_0}$ and $h=\frac{2\pi}{2\ell+1}$, we have for all $j=0,\ldots,2\ell$, that
$$\begin{array}{ccl}
\lambda_j &=&\frac{e^{i(\theta_0+hj)\ell}}{Ne^{i(2\ell+1)\theta_0}}\displaystyle{\sum_{k=1}^{2\ell+1}} q^{\frac{(k-\ell-1)^2}{2}}e^{i(\theta_0+hj)k} \\[10pt]
&=&\frac{1}{N}e^{i\left[ lhj-(l+1)\theta_0 \right]}\left\{ e^{i(\theta_0+hj)(\ell+1)} + \displaystyle{\sum_{k=1}^{\ell}} \left( q^{\frac{(k-\ell-1)^2}{2}}e^{i(\theta_0+hj)k} + q^{\frac{k^2}{2}}e^{i(\theta_0+hj)(\ell+1+k)}\right) \right\} \\[10pt]
&=&\frac{1}{N}\left[ 1 + \displaystyle{\sum_{k=1}^{\ell}} q^{\frac{(k-\ell-1)^2}{2}}e^{i\left[\theta_0(k-\ell-1)+hj(k+\ell)\right]} + q^{\frac{k^2}{2}}e^{i\left[\theta_0k+hj(k+2\ell+1)\right]}    \right] \\[10pt]
&=&\frac{1}{N}\left\{ 1 + \displaystyle{\sum_{k=1}^{\ell}} q^{\frac{k^2}{2}} \left[ e^{i\left( \theta_0k+hjk \right)} + e^{i\left[ -k\theta_0 + hj(2\ell +1)-hjk \right]}  \right] \right\} \\[10pt]
&=&\frac{1}{N}\left\{1+ 2\cdot\displaystyle{\sum_{k=1}^{\ell}}  q^{\frac{k^2}{2}}\cos \left[\left( \theta_0+hj \right)k\right] \right\}.
\end{array}$$
In the case $N=2\ell$, $r=\ell$ and $s=\ell-1$ we have that
\begin{equation}\label{pesosRogersinterpolatexplicitauxiliar}
\begin{array}{ccl}
\lambda_j&=&\frac{z_j^\ell}{N\tau}\displaystyle{\sum_{k=1}^{2\ell}} q^{\frac{(k-\ell)^2}{2}}z_j^k \\[10pt]
&=&\frac{z_j^\ell}{N\tau}\left\{ z_j^\ell + q^{\ell^2/2}z_j^{2\ell}+\displaystyle{\sum_{k=1}^{\ell-1}} q^{k^2/2}\left( z_j^{\ell+k}+z_j^{\ell-k} \right)   \right\} \\[10pt]
&=&\frac{1}{N}\left\{ 1 + q^{\ell^2/2}z_j^{\ell}+2\cdot\displaystyle{\sum_{k=1}^{\ell-1}} q^{k^2/2}\cos\left[(\theta_0+hj)k \right] \right\} \\[10pt]
&=&\frac{1}{N}\left\{1+ (-1)^j q^{\ell^2/2}e^{i\ell\theta_0} + 2\cdot\displaystyle{\sum_{k=1}^{\ell-1}}  q^{\frac{k^2}{2}}\cos \left[\left( \theta_0+hj \right)k\right] \right\},
\end{array}
\end{equation}
so it is clear that the imaginary part of $\lambda_j$ is given by $\frac{(-1)^j}{N}q^{\ell^2/2}\sin(\ell\theta_0)$, that does not depend on $j$, up to the factor $(-1)^j$. In a similar way it can be proved that the case $N=2\ell$, $r=\ell-1$ and $s=\ell$ leads to $\overline{\lambda_j}$, the conjugate of the weight $\lambda_j$ in (\ref{pesosRogersinterpolatexplicitauxiliar}), so (\ref{pesosRogersinterpolatexplicit}) is finally recovered.
\begin{flushright}
$\Box$
\end{flushright}
\end{DT}

In Table \ref{Tab4} we compare the absolute value of the errors obtained when estimating several integrals using the three procedures described above. The values correspond to
\begin{enumerate}
\item Error 1: $\left| I_N^{\omega}(F)-I_{\mu}(F) \right|$ (error 1), with $I_N^{\omega}(F)$ given by (\ref{INRogers}).
\item Error 2: $\left| I_m^{\omega}(F)-I_{\mu}(F) \right|$ (error 2), with $I_m^{\omega}(F)$ given by (\ref{ImRogers}).
\item Error 3: The mixed interpolation-regression method introduced in this paper, that we have computed only when $N-m\geq3$.
\end{enumerate}

In the case of $F_1=F_1(e^{i\theta})$ (first two rows) we can see that the results for |{\em error 1}| ($20$ and $30$-point interpolatory-type rule with uniformly distributed nodes, exact in $\Lambda_{-10,9}$ and $\Lambda_{-15,14}$, respectively) are better that those for |{\em error 2}| ($15$ and $26$-point interpolatory-type rule with nodes not equispaced, but close to the zeros of the para-orthogonal polynomials $B_{15}(z,\tau_{15})$ and $B_{26}(z,\tau_{26})$, respectively, where $\tau_{15} \approx 0.432236577168781 - 0.901760246050699i$ and $\tau_{26} \approx 0.810756105685128 + 0.585384093646459i$, and exact in $\Lambda_{-8,7}$ and $\Lambda_{-13,12}$, respectively). In both cases it holds $r=m+3<N$, so in the mixed interpolation-regression method we approximate the corresponding integral from the direct integration of an approximating Laurent polynomial that belongs to a strict linear subspace of $\Lambda_{-10,9}$ and $\Lambda_{-15,14}$, respectively. For the exponential function, interpolatory-type rules based on equispaced nodes converge to the integral, so we expect worst results when $m$ is far from $N$. This is shown in the third and fourth rows. However, we must highlight the accurate results obtained by using our new procedure. In the third row, the interpolating Laurent polynomials that approximate the $F_1$ in Columns 1-3 belongs in $\Lambda_{-15,14}$, $\Lambda_{-5,4}$ and $\Lambda_{-6,6}$, respectively, whereas in the fourth row, these linear subspaces are $\Lambda_{-10,9}$, $\Lambda_{-4,3}$ and $\Lambda_{-6,5}$, respectively.

In the second example we have considered the following discontinuous function
$$f_2=f_2(\theta)=10\cdot \left(1-2\textrm{heaviside}\left(\theta-\frac{\pi}{2}\right)\right)=\left\{ \begin{array}{rcl} 10 &\textrm{if} & 0 \leq \theta < \frac{\pi}{2}, \\ 0 &\textrm{if} & \theta=\frac{\pi}{2}, \\ -10 &\textrm{if} & \frac{\pi}{2} \leq \theta < 2\pi, \end{array}\right.$$
extended $2\pi$-periodic for $\theta\in\RR$: $$f_2(\theta)=10\cdot \left(1-2\textrm{heaviside}\left(mod(\theta,2\pi)-\frac{\pi}{2}\right)\right).$$ For this particular function, interpolatory-type quadrature formulas based on equidistant nodes do not converge. Some illustrative examples are shown in Table 4,  where the results obtained from the mixed interpolation-regression method clearly improve the two other quadrature rules, for different values of the parameters $\theta_0$, $q$, $m$ and $N$.

We must recall here that the function $L$ solves the minimizing problem (\ref{QN-mProblem}), with respect to the discrete 2-norm $\parallel \cdot \parallel_2$ on $Z_N \backslash \Upsilon_m=\left\{ \hat{z}_s \right\}_{s=1}^{N-m}$. This is not the absolute error we are showing in Table \ref{Tab4}, so we might eventually find examples where the values in column |{\em error 3}| are greater than those in |{\em error 1}| and |{\em error 2}| columns (this is not our case).

\begin{table}[H]
\centering
\resizebox{1.02\textwidth}{!}{
\begin{tabular}{cc|c|c|c|l}
\cline{3-5}
& & $\;$ {\bf |error 1|} $\;$ & $\;${\bf |error 2|}$\;$ & $\;${\bf |error 3|}$\;$ \\ \cline{1-5}
\multicolumn{1}{ |c  }{} &
\multicolumn{1}{ |c| }{$q=0.8$, $m=15$, $N=20$} & 1.199960761886714E-11 & 1.123694361985777E-09 & 4.384092477781099E-12 & \\ \cline{2-5}
\multicolumn{1}{ |c  }{\multirow{2}{*}{$\begin{array}{c} F_1(z)=e^z \\ \left( \theta_0=\frac{5\pi}{6}\right)\end{array}$}} &
\multicolumn{1}{ |c| }{$q=0.7$, $m=26$, $N=30$} & 7.962238595100411E-15 & 2.245086287498004E-12 & 4.592063569987901E-15 & \\ \cline{2-5}
\multicolumn{1}{ |c  }{} &
\multicolumn{1}{ |c| }{$q=0.5$, $m=10$, $N=30$} & 7.514090228864127E-15 & 0.002518296694125 & 0.575429730148636E-05 & \\ \cline{2-5}
\multicolumn{1}{ |c  }{} &
\multicolumn{1}{ |c| }{$q=0.25$, $m=8$, $N=20$} & 5.614802271172991E-15 & 0.018622482512379 & 4.138481385928899E-04 & \\ \cline{1-5}
\multicolumn{1}{ |c  }{} &
\multicolumn{1}{ |c| }{$q=0.01$, $m=10$, $N=24$, $\theta_0=\frac{\pi}{3}$, $r=15$} & 0.496443904040067 & 0.481577391899809 & 0.024444307331453 & \\ \cline{2-5}
\multicolumn{1}{ |c  }{\multirow{2}{*}{$f_2(\theta)$}} &
\multicolumn{1}{ |c| }{$q=0.01$, $m=14$, $N=24$, $\theta_0=\frac{\pi}{3}$, $r=19$} & 0.496443904040067 & 0.334634570232594 & 0.004021991699727 & \\ \cline{2-5}
\multicolumn{1}{ |c  }{} &
\multicolumn{1}{ |c| }{$q=0.05$, $m=40$, $N=60$, $\theta_0=0$, $r=43$} & 0.240735367870585 & 0.728350506052769 & 0.002776303200925 & \\ \cline{2-5}
\multicolumn{1}{ |c  }{} &
\multicolumn{1}{ |c| }{$q=0.5$, $m=10$, $N=24$, $\theta_0=\frac{7\pi}{6}$, $r=14$} & 1.233584727169714 & 0.077374941718200 & 0.013642381379640 & \\ \cline{2-5}
\multicolumn{1}{ |c  }{} &
\multicolumn{1}{ |c| }{$q=0.8$, $m=16$, $N=30$, $\theta_0=0$, $r=19$} & 1.769614754650557 & 0.372801718323300 & 0.047430868386942 & \\ \cline{2-5}
\multicolumn{1}{ |c  }{} &
\multicolumn{1}{ |c| }{$q=0.9$, $m=4$, $N=14$, $\theta_0=\frac{\pi}{3}$, $r=7$} & 1.900770300617066 & 3.041862960006026 & 0.184825151185287 & \\ \cline{1-5}
\end{tabular}
}
\caption{Estimation of integrals on the unit circle.}
\label{Tab4}
\end{table}

We present some final examples where the error in the discrete 2-norm $\parallel \cdot \parallel_2$ on $Z_N \backslash \Upsilon_m=\left\{ \hat{z}_s \right\}_{s=1}^{N-m}$ is computed with the purpose of comparing the advantage when considering Hermite instead of Lagrange interpolation, assuming that the values of some derivatives of the function $F$ are known at the interpolating points.

Consider the function $F(z)=\frac{1}{z-\alpha}$, $\alpha=0.8+0.5i \in \DD$. As indicated in the Introduction, polynomial interpolation based on equally spaced nodes on the unit circle may not converge to $F$ uniformly, and can diverge on points on $\TT$ close to the singularity $\alpha \approx \TT$. Here we are taking advantage of the fact that $P_m$, which belongs to a linear subspace of Laurent polynomials of dimension $m$, can also belong to a linear subspace of higher dimension $r$. The extra conditions are the knowledge of the value of some derivatives of the integrand $F$ at some nodes. In our case we have considered $r-m$ first derivatives distributed at some random nodes of $\Upsilon_m$. Some results are presented in Table \ref{Tab5} taking $\theta_0=\frac{\pi}{6}$. This election of $\theta_0$ is intentional since in this case, a node is fixed at $\TT$ that is very close to the pole of $F$; therefore, poor results can be expected. Actually, these poor results will continue in any case for other values of $\theta_0$ since for sufficiently large values of $N$ we will still have nodes of $Z_N$ located close to the pole $\alpha$. Moreover, observe that as $q$ approaches $1^-$, the zeros of para-orthogonal polynomials tend to accumulate in a region of $\TT$ close to $\alpha$, and therefore, even worse results are seen in this case.

\begin{table}[H]
\centering
\resizebox{0.6\textwidth}{!}{
\begin{tabular}{|c|c|c|c|c|c|}
\hline
$q$ & $N$ & $m$ & $r$ & Error 1 & Error 2 \\ \hline
$0.2$ & $20$ &
$10$
&
$15$ &
$18.655423866409901$ &
$12.770132240717604$
\\ \hline
$0.5$ & $20$ &
$12$
&
$15$ &
$56.859914950230198$ &
$29.809701020432474$
\\ \hline
$0.8$ & $20$ &
$10$
&
$15$ &
$1087.009337208344$ &
$263.9153772620719$
\\ \hline
\end{tabular}
}
\caption{The error in the discrete 2-norm $\parallel \cdot \parallel_2$ on $Z_N \backslash \Upsilon_m=\left\{ \hat{z}_s \right\}_{s=1}^{N-m}$ in the estimation of integrals on the unit circle by using Lagrange ({\em error 1}) and Hermite ({\em error 2}) interpolation. We take $F(z)=\frac{1}{z-\alpha}$ with $\alpha=0.8+0.5i \in \DD$ and $\theta_0=\frac{\pi}{6}$.}
\label{Tab5}
\end{table}

However, convergence is obtained when the pole $\alpha$ is far from $\TT$, see Table \ref{Tab6}.

\begin{table}[H]
\centering
\resizebox{0.6\textwidth}{!}{
\begin{tabular}{|c|c|c|c|c|c|}
\hline
$q$ & $N$ & $m$ & $r$ & Error 1 & Error 2 \\ \hline
$0.2$ & $30$ &
$20$
&
$23$ &
$5.153496042339725E-05$ &
$2.339569540696130E-05$
\\ \hline
$0.2$ & $40$ &
$20$
&
$25$ &
$4.314396224690475E-05$ &
$1.151876855688029E-05$
\\ \hline
$0.5$ & $20$ &
$8$
&
$13$ &
$0.128047686237509$ &
$0.067008452668312$
\\ \hline
$0.7$ & $20$ &
$15$
&
$18$ & $0.003232788251804$ &
$0.002942034709498$
\\ \hline
\end{tabular}
}
\caption{The error in the discrete 2-norm $\parallel \cdot \parallel_2$ on $Z_N \backslash \Upsilon_m=\left\{ \hat{z}_s \right\}_{s=1}^{N-m}$ in the estimation of integrals on the unit circle by using Lagrange ({\em error 1}) and Hermite ({\em error 2}) interpolation. We take $F(z)=\frac{1}{z-\alpha}$ with $\alpha=(1+i)/5 \in \DD$ and $\theta_0=0$.}
\label{Tab6}
\end{table}

The use of Hermite instead of Lagrange interpolation is pointed out in a final example, by taking $F(z)=e^{z/2}$. Results are shown in Table \ref{Tab7}.

%

\begin{table}[H]
\centering
\resizebox{0.6\textwidth}{!}{
\begin{tabular}{|c|c|c|c|c|c|}
\hline
$q$ & $N$ & $m$ & $r$ & Error 1 & Error 2 \\ \hline
$0.2$ & $30$ &
$10$
&
$13$ &
$2.915815134358195E-04$ &
$5.430000289552126E-05$
\\ \hline
$0.5$ & $30$ &
$15$
&
$19$ &
$7.383053777629033E-06$ &
$7.036885020603747E-08$
\\ \hline
$0.9$ & $30$ &
$20$
&
$23$ &
$1.338534424565065E-06$ &
$4.648114157860821E-09$
\\ \hline
\end{tabular}
}
\caption{The error in the discrete 2-norm $\parallel \cdot \parallel_2$ on $Z_N \backslash \Upsilon_m=\left\{ \hat{z}_s \right\}_{s=1}^{N-m}$ in the estimation of integrals on the unit circle by using Lagrange ({\em error 1}) and Hermite ({\em error 2}) interpolation. We take $F(z)=e^{z/2}$ with $\theta_0=0$.}
\label{Tab7}
\end{table}

\section{Concluding remarks}\label{sec5}

As in many practical situations, when you are dealing with the integration  of a complex function $F$ with respect to a positive Borel measure $\mu$ supported on the unit circle it is assumed the knowledge of such a function in  a finite number of points, uniformly distributed on the unit circle. In our contribution we focus our search on an approximating Laurent polynomial $L$ to $F$ by using the Hermite interpolation in a set of these points that mimic the zeros of a para-orthogonal polynomial with respect to $\mu$, and the values of $F$ at the remaining nodes to improve the accuracy of the approximation by a process of simultaneous complex regression.  Numerical experiments illustrate our powerful methods in comparison with the standard Gaussian quadrature formulas on the unit circle (Szeg\H{o} rules) in terms of the estimates of the absolute errors in the above approximations to the integrals.\\

\section{Acknowledgements.}

The work of Lidia Fernández has been supported by grants PID2023-149117NB-I00, PID2024-155133NB-I00 and CEX 2020-001105-M funded by MICIU/AEI/10.13039/501100011033 and ERDF/EU, Spain.

The work of Francisco Marcell\'an  has been supported by the research project PID2024-155133NB-I00, \emph{Ortogonalidad, Aproximaci\'on e Integrabilidad: Aplicaciones en Procesos Estoc\'asticos Cl\'asicos y Cu\'anticos} funded by MCIU/AEI, Spain.


\begin{thebibliography}{00}

\begin{small}

\bibitem[]{Alicia1}
Berriochoa, E. and Cachafeiro, A.,
  \textit{Algorithms for solving Hermite interpolation problems using the Fast Fourier Transform},
  J. Comput. Appl. Math. 235 (2010), no. 4,  882--894.

\bibitem[]{Alicia2}
Berriochoa, E., Cachafeiro, A. and Díaz, J.,
  \textit{Gibbs phenomenon in the Hermite interpolation on the circle},
  Appl. Math. Comp. 253 (2015), 274--286.

\bibitem[]{Alicia3}
Berriochoa, E., Cachafeiro, A. and García Amor, J.M.,
  \textit{Gibbs-Wilbraham oscillation related to an Hermite interpolation problem on the unit circle},
  J. Comput. Appl. Math. 344 (2018), 657--675.

\bibitem[]{BXu}
Boyd, J. P. and Xu, F.,
  \textit{Divergence (Runge Phenomenon) for least-squares polynomial approximation on an equispaced grid and Mock–Chebyshev subset interpolation},
  Appl. Math. Comput. 210 (2009)  no. 1, 158--168.

\bibitem[]{RuyD}
  Cruz-Barroso, R. and Delvaux, S.,
  \textit{Orthogonal Laurent polynomials on the unit circle and snake-shaped matrix factorizations},
  J. Approx. Theory 161 (2009), no. 1, 65--87.

\bibitem[]{CMich}
  Cruz-Barroso, R. and Díaz-Bautista, M.,
  \textit{Szeg{\H{o}} quadrature formulas associated with the weight function $1+\cos(m\theta)$ with application to the approximation of certain Christoffel transformations on the interval},
  Numer. Algorithms,  Published online: 01 May 2025, 40 pages.

\bibitem[]{CDP}
  Cruz-Barroso, R., Díaz Mendoza, C. and Perdomo P{\'{\i}}o, F.,
  \textit{Szeg{\H{o}}-type quadrature formulas},
  J. Math. Anal. Appl. 455 (2017) no. 1,  592--605.

\bibitem[]{RuyRogers}
  Cruz-Barroso, R., González-Vera, P., C. and Perdomo P{\'{\i}}o, F.,
  \textit{Quadrature formulas associated with Rogers-Szeg\H{o} polynomials},
  Comput. Math. Appl. 57 (2009), no. 2,  308--323.

\bibitem[]{Daruis}
 Daruis, L. and González-Vera, P.,
  \textit{A note on Hermite-Fejér interpolation for the unit circle},
  Appl. Math. Lett. 14 (2001), no. 8,  997--1003.

\bibitem[]{Del1}
  De Marchi, D., Dell'Accio, F. and Mazza, M.,
  \textit{On the constrained mock-Chebyshev least squares},
  J. Math. Anal. Appl. 280  (2015) , no. 1, 94--109.

\bibitem[]{Del2}
  Dell'Accio, F., Di Tommaso, F. and Nudo, F.,
  \textit{Constrained mock-Chebyshev least squares quadrature},
  Appl. Math. Letters 134 (2022), Paper no. 108328, 10 pp.

\bibitem[]{Del3}
  Dell'Accio, F., Di Tommaso, F. and Nudo, F.,
  \textit{Generalizations of the constrained mock-Chebyshev least squares in two variables: Tensor product vs total degree polynomial interpolation},
  Appl. Math. Letters 125 (2022), Paper no. 107732, 8 pp.

\bibitem[]{Del4}
  Dell'Accio, F., Mezzanotte, D., Nudo, F. and Occorsio, D.,
  \textit{New results on the constrained mock-Chebyshev least squares operator},
  J. Comput. Appl. Math. 473 (2026), Article 116859, 14 pp.

\bibitem[]{Del5}
  Dell'Accio, F., Mezzanotte, D., Nudo, F. and Occorsio, D.,
  \textit{Product integration rules by constrained mock-Chebyshev least squares operator},
  BIT Numerical Mathematics 63 (2023) no. 2, Paper no. 24, 19  pp.

\bibitem[]{Del6}
  Dell'Accio, F., Marcellán, F. and Nudo, F.,
  \textit{A quadrature formula on triangular domains via an interpolatior-regression approach},
  Appl. Math. Letters 163 (2025), Article 109414, 6 pp.

\bibitem[]{Del7}
  Dell'Accio, F., Marcellán, F. and Nudo, F.,
  \textit{A novel interpolation-regression approach for function approximation on the disk and its application to cubature formulas},
Adv. Comput. Math. 51 (2025), no. 6, Article 53, 27 pp.

\bibitem[]{Del8}
  Dell'Accio, F., Marcellán, F. and Nudo, F.,
  \textit{An extension of a mixed interpolation-regresion method using zeros of orthogonal polynomials},
  J. Comput. Appl. Math. 450 (2024), Article 116010, 10 pp.

\bibitem[]{Del9}
  Dell'Accio, F., Marcellán, F. and Nudo, F.,
  \textit{Constrained mock-Chebyshev least squares approximation for Hermite interpolation},
  Numer. Algorithms.  Published online: 22 March 2025.

\bibitem[]{Gau}
  Gautschi, W.,
  \textit{Orthogonal polynomials: Computation and Approximation},
  Numerical Mathematics and Scientific Computation. Oxford University Press, New York, 2004.

\bibitem[]{Gol}
Golinskii, L.,
\textit{Quadrature formula and zeros of para-orthogonal polynomials on the unit circle},
Acta Math. Hungar. 96 (2002), no. 3, 169--186.

\bibitem[]{JNT}
  Jones, W.B., Nj{\aa}stad, O. and Thron, W.J.,
  \textit{Moment theory, orthogonal polynomials, quadrature, and continued fractions associated with the unit circle},
  Bull. London Math. Soc. 21 (1989),  no. 2,  113--152.

\bibitem[]{Miller}
  Miller, K.S.,
  \textit{Complex Linear Least Squares},
  SIAM Review 15 (1973), no. 4,  706-726.

\bibitem[]{Rei}
  Reichel, L.,
  \textit{On polynomial approximation in the uniform norm by the discrete least squares method},
  BIT 26  (1986), no. 3, 349--368.

\bibitem[]{Santos}
  Nj{\aa}stad, O. and Santos-León, J.C.,
  \textit{Domain of validity of Szeg{\H{o}} quadrature formulas},
  J. Comput. Appl. Math. 202 (2007), no. 2, 440--449.

\bibitem[]{SimonBk}
  Simon, B.,
  \textit{Orthogonal Polynomials on the Unit Circle. Part 1: Classical Theory},
  Amer. Math. Soc. Colloq. Publ.  54, Amer. Math. Soc. Providence, Rhode Island, 2005.

\bibitem[]{Sz}
  Szeg{\H{o}}, G.,
  \textit{Orthogonal Polynomials},
  Amer. Math. Soc.  Colloq. Publ.  23, Amer. Math. Soc. Providence, Rhode Island, 1975. Fourth Edition

\end{small}

\end{thebibliography}
\end{document}